\documentclass[twoside,12pt]{article}

\usepackage{times,amsmath,amssymb,amsthm,cite}

\newcommand{\subj}[1]{\par\noindent{\bf AMS Subject Classifications: }#1.}
\newcommand{\keyw}[1]{\par\noindent{\bf Keywords: }#1.}
\numberwithin{equation}{section}
\numberwithin{figure}{section}

\theoremstyle{definition}

\theoremstyle{remark}

\everymath{\displaystyle}
\date{}
\begin{document}

\title{\center\Large\bf On Second Solutions to Second-Order Difference Equations}
\author{{\bf William C. Parke$^{\dagger}$} and {\bf Leonard C. Maximon\,$^{\dagger\,\ddagger}$}\vspace{0.4mm} \\ 
        $\dagger${\small The George Washington University}\\
        {\small Department of Physics}\\
        {\small Washington DC 20052 USA}\\
        {\tt {\small wparke@gwu.edu}} {\footnotesize and} {\tt {\small max@gwu.edu}} \vspace{-2mm}
\and
        $\ddagger${\small Arizona State University}\\ 
        {\small Department of Physics}\\
        {\small Tempe AZ 85287 USA}
}
\maketitle
\thispagestyle{empty}

\begin{abstract}
We investigate and derive second solutions to linear homogeneous second-order difference equations using a variety of methods, in each case going beyond the purely formal solution and giving explicit expressions for the second solution. We present a new implementation of d'Alembert's reduction of order method, applying it to linear second-order recursion equations. Further, we introduce an iterative method to obtain a general solution, giving two linearly independent polynomial solutions to the recurrence relation. In the case of a particular confluent hypergeometric function for which the standard second solution is not independent of the first, i.e. the solutions are degenerate, we use the corresponding differential equation and apply the extended Cauchy-integral method to find a polynomial second solution for the difference equation. We show that the standard d'Alembert method also generates this polynomial solution. 
\end{abstract}

\subj{33-02, 33C15, 39-02, 39A06}

\vspace{2mm}

\keyw{Confluent hypergeometric, difference equations, differential equations,
recurrence relations, second solutions, polynomial solutions, extended Cauchy-integral method}

\newcommand{\p} {\mkern-6mu+\mkern-4mu}
\newcommand{\m} {\mkern-6mu-\mkern-4mu}

\interfootnotelinepenalty=10000

\newpage

\section{Introduction}

\makeatletter
    \def\tagform@#1{\maketag@@@{\normalsize(#1)\@@italiccorr}}
\makeatother

There are a number of distinct methods for generating a second independent solution to a
second-order linear differential equation when one solution is known. These include
1) The extended Cauchy-integral method, 
2) The application of D'Alembert's reduction of order method, and
3) Recursion using selected starting values.

In this paper we present analogous methods for finding explicit forms for a second independent
solution (also called ``solution of the second kind'') to linear second-order difference equations
when a first solution is known.\footnote{For a general discussion of cases when difference
equations with polynomial coefficients can be solved in terms of polynomials, rational functions, and hypergeometric functions, see the thesis of Christian Weixlbaumer \cite{Weix}.}
Our aim in each case is to go beyond the purely formal solution
and derive explicit expressions for the second solution.

For the confluent hypergeometric function, with the first solution given by\linebreak $\,_1\!F_1(a;b;x)$,
the standard second solution is not independent of the first when the first parameter takes on
negative integer values. In this case we give the explicit polynomial second solution
when the second parameter is a positive integer.
These polynomials also arise in second solutions to the confluent hypergeometric differential equation.

\section{Reduction of order method for the second solution of recurrence relations}

D'Alembert's reduction of order technique is widely used to find second solutions to second-order differential equations. In this section we apply an
analogous technique to the general second-order linear homogeneous recurrence
relation expressed by
\begin{equation}
a_{n}y_{n+2}+b_{n}y_{n+1}+c_{n}y_{n}=0\   \label{gen-rec}
\end{equation}
where $a_{n}c_{n}\ne 0.$

Suppose we know one solution of the above recurrence relation to be
\begin{equation}
y_{n}^{( 1) }=f_{n} \ .
\end{equation}
Then an ansatz for a second solution, possibly independent of the first, is
\begin{equation}
y_{n}^{( 2) }=w_{n}f_{n} \ .
\end{equation}
The difference equation for $w_{n}$ can be solved
explicitly.

Assuming $y_{n}^{( 2) }$\ to be a solution of the recurrence
relation, and using $\Delta w_n\equiv w_{n+1}-w_n,$ the $w_{n}$ will satisfy
\begin{equation}
a_{n}\left( w_{n}+2\Delta w_{n}+\Delta ^{2}w_{n}\right) f_{n+2}+b_{n}\left(
w_{n}+\Delta w_{n}\right) f_{n+1}+c_{n}w_{n}f_{n}=0
\end{equation}
which becomes, since $f_{n}$ is a solution,
\begin{equation}
a_{n}\left( 2\Delta w_{n}+\Delta ^{2}w_{n}\right) f_{n+2}+b_{n}\left( \Delta
w_{n}\right) f_{n+1}=0 \ .
\end{equation}
Let
\begin{equation}
u_{n}=\Delta w_{n} \ .
\end{equation}
Then
\begin{eqnarray}
a_{n}f_{n+2}\left( 2u_{n}+\Delta u_{n}\right) +b_{n}\left( u_{n}\right)
f_{n+1} &=&0 \nonumber\\
a_{n}f_{n+2}\left( u_{n}+u_{n+1}\right) +b_{n}f_{n+1}u_{n} &=&0\ .
\end{eqnarray}
The $u_{n}$ satisfy a first-order difference equation (assuming none of the $%
f_{n}$ vanish):
\begin{eqnarray}
u_{n+1} &=&-\left( 1+\frac{b_{n}f_{n+1}}{a_{n}f_{n+2}}\right) u_{n} \\
&=&\left( \frac{c_{n}f_{n}}{a_{n}f_{n+2}}\right) u_{n}\ .
\end{eqnarray}
Iterating,
\begin{equation}
u_{n}=\prod_{l=0}^{n-1}\left( \frac{c_{l}f_{l}}{a_{l}f_{l+2}}\right) u_{0}
\end{equation}
or
\begin{equation}
\Delta w_{k}=\frac{f_{0}f_{1}}{f_{k}f_{k+1}}\prod_{l=0}^{k-1}\left( \frac{%
c_{l}}{a_{l}}\right) u_{0}\ .
\end{equation}
Summing the $\Delta w_{k}$\ gives
\begin{equation}
w_{n}=w_{0}+u_{0}\sum_{k=0}^{n-1}\frac{f_{0}f_{1}}{f_{k}f_{k+1}}%
\prod_{l=0}^{k-1}\left( \frac{c_{l}}{a_{l}}\right)\ .
\end{equation}
(Following convention, we
take products that have an upper limit smaller than the lower limit to be unity and
sums that have an upper limit smaller than the lower limit as vanishing.)
We can drop $w_{0}$ and the term $k=0$ in the sum, as they will reproduce
the first solution, and also select the overall factor to be independent of the
index $n$, since the general solution is constructed by a sum of an
arbitrary constant times each of two independent solutions. We write our second solution as
\begin{equation}
y_{n}^{( 2) }=f_{1}f_{n}\sum_{k=1}^{n-1}\frac{1}{f_{k}f_{k+1}}%
\prod_{l=0}^{k-1}\left( \frac{c_{l}}{a_{l}}\right) \ . \label{2nd-sol1}
\end{equation}

In the next section we evaluate the second solution as given in (\ref{2nd-sol1}) for a few simple
examples. We then consider the reduction of the sum in (\ref{2nd-sol1}), and derive a recurrence relation
for the second solution.

\subsection{Examples of second solutions}

In this section, we evaluate the solution Eq. (\ref{2nd-sol1}) for some simple examples.

\begin{enumerate}
\item  First, the case with constant coefficients,
\begin{equation}
ay_{n+2}+by_{n+1}+cy_{n}=0
\end{equation}
and a double root to the characteristic equation $ar^2+br+c=0$, i.e., $%
r=-b/2a $ with $b^{2}=4ac.$ One solution to the given recurrence relation is
\begin{equation}
y_{n}^{( 1) }=r^{n} \ .
\end{equation}
Our second solution is then
\begin{eqnarray}
y_{n}^{( 2) } &=&f_{1}f_{n}\sum_{k=1}^{n-1}\frac{1}{f_{k}f_{k+1}}%
\prod_{l=0}^{k-1}\left( \frac{c_{l}}{a_{l}}\right)  \nonumber\\
&=&r^{n+1}\sum_{k=1}^{n-1}\frac{1}{r^{2k+1}}\prod_{l=0}^{k-1}\left(
r^{2}\right) \nonumber \\
&=&r^{n+1}\sum_{k=1}^{n-1}\frac{1}{r^{2k+1}}r^{2k}=r^{n}\sum_{k=1}^{n-1}1=%
\left( n-1\right) r^{n} \ .
\end{eqnarray}
The expression contains a linear combination of the first and a second
solution, which we can take as $nr^n$.

Now, suppose we do not have a double root, but rather
\begin{equation}
r_{\pm }=-b/2a\pm \sqrt{\left( b/2a\right) ^{2}-c/a}\ .
\end{equation}
Starting with
\begin{equation}
y_{n}^{( 1) }=r_{+}^{n} \ ,
\end{equation}
we have a second solution expressed as
\begin{eqnarray}
y_{n}^{( 2) } &=&r_{+}^{n+1}\sum_{k=1}^{n-1}\frac{1}{r_{+}^{2k+1}}%
\prod_{l=0}^{k-1}\left( \frac{c}{a}\right) \nonumber \\
&=&r_{+}^{n}\sum_{k=1}^{n-1}\frac{1}{r_{+}^{2k}}\left( \frac{c}{a}\right)
^{k} \nonumber \\
&=&r_{+}^{n}\sum_{k=1}^{n-1}\left( \frac{c}{a}\frac{1}{r_{+}^{2}}\right) ^{k}\nonumber \\
&=&r_{+}^{n}\frac{1}{\frac{c}{a}\frac{1}{r_{+}^{2}}-1}\left( \left( \frac{c}{%
a}\frac{1}{r_{+}^{2}}\right) ^{n}-1\right) -r_{+}^{n} \nonumber \\
&=&\frac{ar_{+}^{2}}{c-ar_{+}^{2}}\left( \left( \frac{c}{a}\frac{1}{r_{+}}%
\right) ^{n}-r_{+}^{n}\right) -r_{+}^{n} \ .
\end{eqnarray}
The above is simplified with
\begin{eqnarray}
c-ar_{+}^{2} &=&br_{+}+2c \nonumber \\
&=&-b^{2}/2a+b\sqrt{\left( b/2a\right) ^{2}-c/a}+2c \ ,
\end{eqnarray}
and
\begin{eqnarray}
r_{+}+r_{-} &=&-b/a \ ,  \\
r_{+}r_{-} &=&\left( -b/2a+\sqrt{\left( b/2a\right) ^{2}-c/a}\right) \left(
-b/2a-\sqrt{\left( b/2a\right) ^{2}-c/a}\right) \ , \nonumber  \\
&=&c/a \ , \\
ar_{+}/c &=&\frac{1}{r_{-}} \ .
\end{eqnarray}
So, the `new part' of our solution becomes
\begin{eqnarray}
y_{n}^{\left( 2^{\prime }\right) } &=&\left( \frac{c}{a}\frac{1}{r_{+}}%
\right) ^{n} \nonumber \\
&=&r_{-}^{n} \ .
\end{eqnarray}
We recover the second solution from the first solution.

\item  As an example with a factorial solution, consider
\begin{equation}
y_{n+2}-\left( n+1\right) y_{n+1}-\left( n+1\right) y_{n}=0 \ .
\end{equation}
One solution is
\begin{equation}
y_{n}^{\left( 1\right) }=n! \ .
\end{equation}
Our second solution will be
\begin{eqnarray}
y_{n}^{( 2) } &=&f_{1}f_{n}\sum_{k=1}^{n-1}\frac{1}{f_{k}f_{k+1}}%
\prod_{l=0}^{k-1}\left( \frac{c_{l}}{a_{l}}\right) \ ,  \\
&=&n!\sum_{k=1}^{n-1}\frac{1}{k!\left( k+1\right) !}\prod_{l=0}^{k-1}\left(
-1\right)\left( l+1\right) \ ,  \\
&=&n!\sum_{k=1}^{n-1}\frac{\left( -1\right) ^{k}}{\left( k+1\right) !}%
=n!\sum_{l=2}^{n}\frac{\left( -1\right) ^{l-1}}{l!} \ .
\end{eqnarray}
A second solution independent of the first can be taken as
\begin{equation}
y_{n}^{( 2^{\prime }) }=n!\sum_{l=0}^{n}\frac{\left( -1\right)
^{l}}{l!}\ .
\end{equation}

\item  An example giving a harmonic-number solution comes from solving
\begin{equation}
\left( n+2\right) y_{n+2}-\left( 2n+3\right) y_{n+1}+\left( n+1\right)
y_{n}=0
\end{equation}
which is the same as the difference relation
\begin{equation}
\left( n+2\right) \Delta ^{2}y_{n}+\Delta y_{n}=0
\end{equation}
(where $\Delta y_n\equiv y_{n+1}-y_n$). Evidently, one solution is just a
constant.

Starting with the constant solution: $f_k=1$, we write our second solution
\begin{eqnarray}
y_{n}^{( 2) } &=&f_{1}f_{n}\sum_{k=1}^{n-1}\frac{1}{f_{k}f_{k+1}}%
\prod_{l=0}^{k-1}\left( \frac{c_{l}}{a_{l}}\right) \nonumber \\
&=&\sum_{k=1}^{n-1}\prod_{l=0}^{k-1}\left( \frac{l+1}{l+2}\right) \nonumber \\
&=&\sum_{k=1}^{n-1}\frac{1}{k+1}=\sum_{l=1}^{n}\frac{1}{l}-1
\end{eqnarray}
so a new (second) solution is
\begin{equation}
y_{n}^{( 2^{\prime }) }=\sum_{l=1}^{n}\frac{1}{l}
\end{equation}
which is the harmonic number $H_n$ as promised.

\item  Another less trivial example comes from
\begin{equation}
\left( n+1\right) y_{n+2}-\left( n^{2}+7n+8\right) y_{n+1}+2\left(
n+2\right) \left( n+3\right) y_{n}=0
\end{equation}
which has a solution
\begin{equation}
y_{n}^{( 1) }=2^{n} \ .
\end{equation}
We construct the second solution using
\begin{eqnarray}
y_{n}^{( 2) } &=&f_{1}f_{n}\sum_{k=1}^{n-1}\frac{1}{f_{k}f_{k+1}}%
\prod_{l=0}^{k-1}\left( \frac{c_{l}}{a_{l}}\right) \nonumber \\
&=&2^{n+1}\sum_{k=1}^{n-1}\frac{1}{2^{2k+1}}\prod_{l=0}^{k-1}\frac{2\left(
l+2\right) \left( l+3\right) }{\left( l+1\right) } \nonumber \\
&=&2^{n+1}\sum_{k=1}^{n-1}\frac{2^{k-1}}{2^{2k+1}}\left( k+1\right) \left(
k+2\right) ! \nonumber \\
&=&2^{n-1}\sum_{k=1}^{n-1}\frac{1}{2^{k}}\left( \left( k+3\right) -2\right)
\left( k+2\right) ! \nonumber \\
&=&2^{n-1}\left( \sum_{k=1}^{n-1}\frac{1}{2^{k}}\left( k+3\right)
!-\sum_{k=1}^{n-1}\frac{1}{2^{k-1}}\left( k+2\right) !\right)
\end{eqnarray}
In the second sum, let $k=l+1.$ \ Then
\begin{eqnarray}
y_{n}^{( 2) } &=&2^{n-1}\left( \sum_{k=1}^{n-1}\frac{1}{2^{k}}%
\left( k+3\right) !-\sum_{l=0}^{n-2}\frac{1}{2^{l}}\left( l+3\right) !\right)\nonumber
\\
&=&\left( n+2\right) !-3\cdot 2^{n} \ .
\end{eqnarray}
We see that a second independent solution is
\begin{equation}
y_{n}^{( 2^{\prime }) }=\left( n+2\right) ! \ .
\end{equation}
\end{enumerate}

\subsection{Reduction of the second solution sum}

We have found that the second solution to the linear homogeneous second-order recurrence
relation, Eq.\,(\ref{gen-rec}), can be usefully expressed
as
\begin{equation}
y_{n}^{(2)}=\sum_{k=1}^{n-1}\frac{f_{1}f_{n}}{f_{k}f_{k+1}}\prod_{l=0}^{k-1}%
\frac{c_{l}}{a_{l}} \ .
\end{equation}
This sum may look intimidating, especially when the first
solution is not a simple function of the index or auxiliary parameters
within the coefficients. However, in these cases, by a sequential `peeling
back' on the summation terms, starting with $k=n-2$ and $k=n-1$, the sum can
be made simpler in form. Consider
\begin{eqnarray}
y_{n}^{(2)} &=&\sum_{k=1}^{n-3}\frac{f_{1}f_{n}}{f_{k}f_{k+1}}%
\prod_{l=0}^{k-1}\frac{c_{l}}{a_{l}}+f_{1}f_{n}\left( \frac{1}{f_{n-2}f_{n-1}}+\frac{1}{f_{n-1}f_{n}}\frac{c_{n-2}}{a_{n-2}}\right)\prod_{l=0}^{n-3}\frac{c_{l}}{%
a_{l}} \nonumber \\
&=&\sum_{k=1}^{n-3}\frac{f_{1}f_{n}}{f_{k}f_{k+1}}\prod_{l=0}^{k-1}\frac{%
c_{l}}{a_{l}}+\frac{f_{1}f_{n}}{f_{n-2}f_{n-1}f_{n}}\frac{1}{a_{n-2}}%
\left( a_{n-2}f_{n}+c_{n-2}f_{n-2}\right)\prod_{l=0}^{n-3}\frac{c_{l}}{a_{l}}\ .\nonumber\\
\end{eqnarray}
Now apply
\begin{equation}
a_{n-2}f_{n}+c_{n-2}f_{n-2}=-b_{n-2}f_{n-1}
\end{equation}
to get
\begin{eqnarray}
y_{n}^{(2)} &=&\sum_{k=1}^{n-3}\frac{f_{1}f_{n}}{f_{k}f_{k+1}}%
\prod_{l=0}^{k-1}\frac{c_{l}}{a_{l}}-\frac{f_{1}}{f_{n-2}}\frac{b_{n-2}}{%
a_{n-2}}\prod_{l=0}^{n-3}\frac{c_{l}}{a_{l}} \nonumber \\
&=&\sum_{k=1}^{n-4}\frac{f_{1}f_{n}}{f_{k}f_{k+1}}\prod_{l=0}^{k-1}\frac{%
c_{l}}{a_{l}}+\left( \frac{f_{1}f_{n}}{%
f_{n-3}f_{n-2}}-\frac{f_{1}}{f_{n-2}}\frac{b_{n-2}}{a_{n-2}}\frac{c_{n-3}}{%
a_{n-3}}\right)\prod_{l=0}^{n-4}\frac{c_{l}}{a_{l}} \nonumber \\
&=&\sum_{k=1}^{n-4}\frac{f_{1}f_{n}}{f_{k}f_{k+1}}\prod_{l=0}^{k-1}\frac{%
c_{l}}{a_{l}}+\frac{f_{1}}{f_{n-3}}\frac{1}{a_{n-3}a_{n-2}}\left(
-a_{n-3}c_{n-2}+b_{n-3}b_{n-2}\right) \prod_{l=0}^{n-4}\frac{c_{l}}{a_{l}}\nonumber \\
\end{eqnarray}
and so forth, until the first sum drops to zero terms. The
$f_1$ factor in the second term will then cancel with an $f_1$ in the
denominator, leaving no more denominator factors of $f_k$.

The general form of the result will be
\begin{equation}
y_{n}^{(2)}=\left( \prod_{l=0}^{n-2}\frac{1}{a_{l}}\right) Y_{n} \   \label{P-form}
\end{equation}
where $Y_n$,
satisfying
\begin{equation}
Y_{n+2}+b_{n}Y_{n+1}
+a_{n-1}c_{n}Y_{n} =0 \ ,
\end{equation}
can be expressed in terms the initial values $Y_0$ and $Y_1$ times polynomials in the set of coefficients
$\{a_0,a_1,\cdots;b_0,b_1,\cdots;c_0,c_1,\cdots\}$. These polynomials will be described in section \ref{rules}.

\section{Iterative derivation of the second solutions}

The second solution given by Eq.\,(\ref{P-form}) can also be deduced by direct
iteration. Start with the linear second-order difference equation (\ref{gen-rec}%
). Make the substitution
\begin{equation}
y_{n}=\left(\prod_{l=0}^{n-2}\frac{1}{a_{l}}\right) Y_{n} \ .
\end{equation}
Let $\beta _{n}=-b_{n}$ and $\gamma _{n}=-a_{n-1}c_{n},$ and define $a_{-1}=1
$. Then
\begin{equation}
Y_{n+2}=\beta _{n}Y_{n+1}+\gamma _{n}Y_{n} \ .
\end{equation}

Realizing that the iterated solution will depend on the pair of initial
values, say $Y_{0}$ and $Y_{1}$, we write the second-order difference equation
as a 2x2 matrix equation:
\begin{equation}
\left(
\begin{array}{c}
Y_{n+2} \\
Y_{n+1}
\end{array}
\right) =\left(
\begin{array}{cc}
\beta _{n} & \gamma _{n} \\
1 & 0
\end{array}
\right) \left(
\begin{array}{c}
Y_{n+1} \\
Y_{n}
\end{array}
\right) \ .
\end{equation}
Iteration gives
\begin{equation}
\left(
\begin{array}{c}
Y_{n+2} \\
Y_{n+1}
\end{array}
\right) =\prod_{l=0}^{n}\left(
\begin{array}{cc}
\beta _{l} & \gamma _{l} \\
1 & 0
\end{array}
\right) \left(
\begin{array}{c}
Y_{1} \\
Y_{0}
\end{array}
\right) \ .  \label{Y-iterated}
\end{equation}

\subsection*{Constant coefficient example}

As a special case, when the recurrence coefficients $\beta$ and $\gamma$ do
not depend on their index, we will have
\begin{eqnarray}
Y_{n}&=& Y_{1}\sum_{k=0}^{\left\lfloor \left( n-1\right) /2\right\rfloor }%
\left(\begin{array}{c}n-1-k\\ k\end{array}\right)\beta ^{n-2k-1}\gamma ^{k} \nonumber \\
&&+\ Y_{0}\sum_{k=0}^{\left\lfloor \left( n-1\right) /2\right\rfloor-1 }
\left(\begin{array}{c}n-2-k\\ k\end{array}\right)
\beta ^{n-2k-2}\gamma ^{k+1} \\
Y_{2} &=&Y_{1}\beta +Y_{0}\gamma \nonumber \\
Y_{3} &=&Y_{1}\left( \beta ^{2}+\gamma \right) +Y_{0}\beta \gamma \nonumber \\
Y_{4} &=&Y_{1}\left( \beta ^{2}+2\gamma \right) \beta +Y_{0}\left( \beta
^{2}+\gamma \right) \gamma \nonumber \\
Y_{5} &=&Y_{1}\left( \beta ^{4}+3\beta ^{2}\gamma +\gamma ^{2}\right)
+Y_{0}\left( \beta ^{3}\gamma +2\beta \gamma ^{2}\right) \ . \nonumber
\end{eqnarray}
These solutions must correspond to the simpler looking ones commonly found
from the characteristic equation, namely
\begin{equation}
Y_{n}=c_{1}r_{1}^{n}+c_{2}r_{2}^{n}  \label{Y-roots}
\end{equation}
where
\begin{eqnarray}
r_{1} &=&\frac{\beta }{2}+\sqrt{\frac{\beta ^{2}}{4}+\gamma } \\
r_{2} &=&\frac{\beta }{2}-\sqrt{\frac{\beta ^{2}}{4}+\gamma } \ .
\end{eqnarray}
Expanding the $n^{th}$ power of these roots into a binomial series
gives
\begin{eqnarray}
r_{2}^{n} &=&\left( \frac{\beta }{2}-\sqrt{\frac{\beta ^{2}}{4}+\gamma }%
\right) ^{n}\nonumber \\
&=&\sum_{k=0}^{n}\left( -1\right) ^{k}
\left(\begin{array}{c}n\\ k\end{array}\right)
\left( \frac{\beta }{2}%
\right) ^{n-k}\left( \frac{\beta ^{2}}{4}+\gamma \right) ^{k/2}
\end{eqnarray}
and similarly for $r_{1}^{n}$ but without the $(-1)^k$ factor. Now
\begin{eqnarray}
Y_{2} &=&c_{1}\left( \frac{1}{2}\beta ^{2}+\beta \sqrt{\frac{1}{4}\beta
^{2}+\gamma }+\gamma \right) +c_{2}\left( \frac{1}{2}\beta ^{2}-\beta \sqrt{%
\frac{1}{4}\beta ^{2}+\gamma }+\gamma \right) \nonumber \\
&=&\frac{1}{2}\left( c_{1}+c_{2}\right) \left( \beta ^{2}+2\gamma \right) +%
\frac{1}{2}\left( c_{1}-c_{2}\right) \beta \sqrt{\beta ^{2}+4\gamma } \nonumber \\
&=&\beta Y_{1}+\gamma Y_{0}
\end{eqnarray}
\begin{eqnarray}
Y_{3} &=&c_{1}r_{1}^{3}+c_{2}r_{2}^{3} \nonumber \\
&=&\frac{1}{2}\beta \left( c_{1}+c_{2}\right) \left( \beta ^{2}+3\gamma
\right) +\frac{1}{2}\sqrt{\beta ^{2}+4\gamma }\left( c_{1}-c_{2}\right)
\left( \beta ^{2}+ \gamma \right) \nonumber \\
&=&Y_{1}\left( \beta ^{2}+\gamma \right) +Y_{0}\beta \gamma
\end{eqnarray}
etc.

Eq.\,(\ref{Y-roots}) can be expressed as a matrix equation:
\begin{equation}
\left(
\begin{array}{c}
Y_{n} \\
Y_{n-1}
\end{array}
\right) =\left(
\begin{array}{cc}
r_{1}^{n} & r_{2}^{n} \\
r_{1}^{n-1} & r_{2}^{n-1}
\end{array}
\right) \left(
\begin{array}{c}
c_{1} \\
c_{2}
\end{array}
\right)  \label{Y-c}
\end{equation}
while our solution is expressed as
\begin{equation}\label{Y-beta}
\left(
\begin{array}{c}
Y_{n} \\
Y_{n-1}
\end{array}
\right) =\left(
\begin{array}{cc}
\beta & \gamma \\
1 & 0
\end{array}
\right) ^{n-1}\left(
\begin{array}{c}
Y_{1} \\
Y_{0}
\end{array}
\right)\,\,.
\end{equation}

Our solution Eq.\,(\ref{Y-beta}) can be transformed to the expression Eq.\,(%
\ref{Y-c}) by finding the eigenvalues of the matrix that appears in Eq.\,(%
\ref{Y-beta}). Suppose the matrix $S$ has the property that
\begin{equation}
S\left(
\begin{array}{cc}
\beta & \gamma \\
1 & 0
\end{array}
\right) S^{-1} =\left(
\begin{array}{cc}
\lambda _{1} & 0 \\
0 & \lambda _{2}
\end{array}
\right) \ .
\end{equation}
It is easy to show that one such matrix is
\begin{equation}
S=\left(
\begin{array}{cc}
r_{1}/\gamma & 1 \\
r_{2}/\gamma & 1
\end{array}
\right)
\end{equation}
with eigenvalues $\lambda_1=r_1$ and $\lambda_2=r_2$, so that
\begin{eqnarray}
S\left(
\begin{array}{c}
Y_{n} \\
Y_{n-1}
\end{array}
\right) &=&\left( S\left(
\begin{array}{cc}
\beta & \gamma \\
1 & 0
\end{array}
\right) S^{-1}\right) ^{n-1}S\left(
\begin{array}{c}
Y_{1} \\
Y_{0}
\end{array}
\right) \\
\left(
\begin{array}{cc}
\lambda _{1}/\gamma & 1 \\
\lambda _{2}/\gamma & 1
\end{array}
\right) \left(
\begin{array}{c}
Y_{n} \\
Y_{n-1}
\end{array}
\right) &=&\left(
\begin{array}{cc}
r_{1}^{n-1} & 0 \\
0 & r_{2}^{n-1}
\end{array}
\right) \left(
\begin{array}{cc}
r_{1}/\gamma & 1 \\
r_{2}/\gamma & 1
\end{array}
\right) \left(
\begin{array}{c}
Y_{1} \\
Y_{0}
\end{array}
\right)
\end{eqnarray}
resulting in
\begin{equation}
Y_{n} =\frac{\left( r_{1}^{n}-r_{2}^{n}\right) }{\left( r_{1}-r_{2}\right) }%
Y_{1}+\gamma \frac{\left( r_{1}^{n-1}-r_{2}^{n-1}\right) }{\left(
r_{1}-r_{2}\right) }Y_{0} \ .
\end{equation}
This solution from the recurrence relations is, as expected, polynomial in
the parameters $\beta$ and $\gamma$ (no square-roots!). We can see this by
observing that
\begin{equation}
\frac{r_{1}^{n}-r_{2}^{n}}{r_{1}-r_{2}} =\frac{1}{2^{n-1}}\sum_{l}
\left(\begin{array}{c}n\\2l+1\end{array}\right)
\beta ^{n-2l-1}\left( \beta ^{2}+4\gamma \right) ^{l} \ .
\end{equation}
Note also that as $r_{1}$ approaches $r_{2}$, the ratio above becomes a
derivative, giving solutions proportional to $nr^{n-1}.$

\subsection{Independence of second solution from the first}

The functions $f_{1}$ and $f_{2}$ are linearly dependent if there exists a
relation
\begin{equation}
c_{1}f_{1}\left( n\right) +c_{2}f_{2}\left( n\right) =0
\end{equation}
with $n$ in a defined range, and the constants $c_1$
and $c_2$ are not zero.

The Casorati determinant for the pair $f_{1}\left( n\right) $, $f_{2}\left(
n\right) $ is defined to be
\begin{equation}
\mathcal{C}\left( n+1\right) =\det \left(
\begin{array}{cc}
f_{1}\left( n\right) & f_{2}\left( n\right) \\
f_{1}\left( n+1\right) & f_{2}\left( n+1\right)
\end{array}
\right) \ .
\end{equation}
The $f_{1}$ and $f_{2}$ will be linearly dependent iff $C\left( n\right) =0$
over the range of $n$.

For our solutions,
\begin{equation}
\left(
\begin{array}{c}
Y_{2} \\
Y_{1}
\end{array}
\right) =\left(
\begin{array}{cc}
\beta _{0} & \gamma _{0} \\
1 & 0
\end{array}
\right) \left(
\begin{array}{c}
Y_{1} \\
Y_{0}
\end{array}
\right)
\end{equation}
i.e.,
\begin{eqnarray}
\left(
\begin{array}{c}
f_{1}\left( 1\right) \\
f_{2}\left( 1\right)
\end{array}
\right) &=&\left(
\begin{array}{c}
\gamma _{0} \\
\beta _{0}
\end{array}
\right)   \\
\left(
\begin{array}{c}
f_{1}\left( 2\right) \\
f_{2}\left( 2\right)
\end{array}
\right) &=&\left(
\begin{array}{c}
\beta _{1}\gamma _{0} \\
\beta _{1}\beta _{0}+\gamma _{1}
\end{array}
\right)
\end{eqnarray}
\begin{eqnarray}
\mathcal{C}\left( 2\right) &=&\det \left(
\begin{array}{cc}
f_{1}\left( 1\right) & f_{2}\left( 1\right) \\
f_{1}\left( 2\right) & f_{2}\left( 2\right)
\end{array}
\right)  \nonumber \\
&=&\det \left(
\begin{array}{cc}
\gamma _{0} & \beta _{0} \\
\beta _{1}\gamma _{0} & \beta _{1}\beta _{0}+\gamma _{1}
\end{array}
\right)  \nonumber \\
&=&\gamma _{0}\gamma _{1} \ .
\end{eqnarray}
Similarly,
\begin{eqnarray}
\det \left(
\begin{array}{cc}
f_{1}\left( 2\right) & f_{2}\left( 2\right) \\
f_{1}\left( 3\right) & f_{2}\left( 3\right)
\end{array}
\right)
&=&\det \left( \left(
\begin{array}{cc}
0 & 1 \\
\gamma _{2} & \beta _{2}
\end{array}
\right) \left(
\begin{array}{cc}
0 & 1 \\
\gamma _{1} & \beta _{1}
\end{array}
\right) \left(
\begin{array}{cc}
0 & 1 \\
\gamma _{0} & \beta _{0}
\end{array}
\right) \right)  \nonumber \\
&=&-\gamma _{0}\gamma _{1}\gamma _{2} \ .
\end{eqnarray}
More generally, we will have
\begin{eqnarray}
\mathcal{C}\left( n+1\right) &=&\det \left(
\begin{array}{cc}
f_{1}\left( n\right) & f_{2}\left( n\right) \\
f_{1}\left( n+1\right) & f_{2}\left( n+1\right)
\end{array}
\right) \nonumber \\
&=&\left( -1\right) ^{n+1}\prod_{l=0}^{n}\gamma _{l} \ .
\end{eqnarray}
The solutions up to $y_{n+1}$ will be independent as long as $\gamma_k\ne0$
for $1\le k \le n$.

\subsection{Rules for constructing the general solutions}\label{rules}

In view of $\beta_n=-b_n$ and $\gamma_n=-a_{n-1}c_n$, and in order to
simplify keeping track of indices, we define
\begin{equation*}
\gamma_{n-1,n}\equiv \gamma_n \ .
\end{equation*}
Because we let $a_{-1}=1,$ we have $\gamma_{-1,0}=-c_0.$

Eq.\,(\ref{Y-iterated}) becomes
\begin{equation}
\left(
\begin{array}{c}
Y_{n+2} \\
Y_{n+1}
\end{array}
\right) =\prod_{k=0}^{n}\left(
\begin{array}{cc}
\beta _{k} & \gamma _{k-1,k} \\
1 & 0
\end{array}
\right) \left(
\begin{array}{c}
Y_{1} \\
Y_{0}
\end{array}
\right) \ .  \label{Y-it}
\end{equation}

By examining the solution (\ref{Y-it}) expanded into polynomials in $\beta_n$
and $\gamma_{n-1,n}$, the following rules for constructing $Y_{n+2}$ in
terms of the initial values $Y_{1}$ and $Y_{0}$ apply:
\begin{enumerate}
\item In the iterated solution of $Y_{n+1}=\beta _{n}Y_{n}+\gamma%
_{n-1,n}Y_{n-1}$, expressed as
\begin{equation}
Y_{n+1}=P_{n}Y_{1}+\gamma _{-1,0}Q_{n-1}Y_{0}\ ,
\end{equation}
the factors $P_n$ and $Q_n$ will be polynomials in the coefficients $%
\beta _{l}$ and $\gamma _{l-1,l},$ homogeneous of degree $n$
in the sense that under the scaling $\beta_l \rightarrow \lambda\beta_l$ and
$\gamma_{l-1,l}\rightarrow \lambda^2 \gamma_{l-1,l}$, we will have $P_n\rightarrow \lambda^n P_n$
and $Q_n\rightarrow \lambda^n Q_n.$
\item There will be a Fibonacci number $F_{n+1}$
of terms in the polynomial $%
P_{n}$. (This can be seen by substituting ones for $\beta_n$ and $\gamma_{n-1,n}$
in Eq.\,(\ref{Y-it}). Here, $F_n=(( (1+\sqrt{5})/2)^n-(1-\sqrt{5})/2)^n)/\sqrt{5}\rightarrow
\{1,1,2,3,5,8,13,\cdots\}.$ )
\item The polynomial $P_{n}$ is constructed as follows:
\begin{enumerate}
\item\label{r3a} For $n$ even, form an even number of initially unindexed factors of $%
\beta $, starting with $n$ such factors (with no $\gamma $ factor), and then
work down to zero $\beta$ factors. For each term with $2k$ factors of $\beta $,
put in $ n /2-k$ factors of $\gamma $, at first unindexed.
Enumerate the terms having a given number of $\beta $ factors to produce all
possible positions of the $\beta ^{\prime }s$\ among the $\gamma $ factors.
The last term with no $\beta $ factors will have $ n /2$
factors of $\gamma $. Now put indices on the $\beta $ and $\gamma $ factors,
sequentially, from $0$ to $n-1.$
\item\label{r3b} For $n$ odd, form an odd number of initially unindexed factors of $%
\beta $ in each term, starting with $n$ factors (with no $\gamma $
factor), and working down to one $\beta$ factor. For each term with $2k+1$ factors
of $\beta $, put $(n-1)/2-k$ factors of $\gamma ,$ at first unindexed. Enumerate
the terms having a given number of $\beta $\ factors to produce all possible
positions of the $\beta ^{\prime }s$ among the $\gamma $ factors. The last
set of terms will have just one $\beta $ factor. Now put indices on the $%
\beta $ and $\gamma $ factors, sequentially, from $0$ to $n-1.$
\end{enumerate}
\item\label{r4} The polynomial $Q_{n}$ is constructed just like $P_{n}$, except that the
indices run from $1$ to $n$ instead of from $0$ to $n-1.$
\end{enumerate}

To exercise these rules, let's write out an example for $n=3$ to get $Y_5$
from Eq.\,(\ref{Y-iterated}).  The functions
$P_4$ and $Q_3$ will be polynomials in the $\beta$ and $\gamma$ of
degree no higher than $4$ and $3$, respectively. There will be $F_5=5$ terms
in $P_4$, and $F_4=3$ terms in $Q_3$. Now, from rule \ref{r3a}, we start
constructing $P_4$ by writing the set
\begin{equation}
\{\beta \beta \beta \beta, \beta \beta \gamma, \beta \gamma \beta, \gamma \beta\beta,
\gamma\gamma\} \ .
\end{equation}
Now decorate sequentially with indices and add:
\begin{equation}
\beta_0 \beta_1 \beta_2 \beta_3+ \beta_0 \beta_1 \gamma_{2,3}+ \beta_0 \gamma_{1,2} \beta_3+ \gamma_{0,1} \beta_2\beta_3+\gamma_{0,1}\gamma_{2,3} \ .
\end{equation}
This is $P_4$.  For $Q_3$, we apply rule \ref{r3b} to generate the set
\begin{equation}
\{\beta \beta \beta,  \beta \gamma, \gamma \beta\} \ .
\end{equation}
With indices according to rule \ref{r4}, the set produces
\begin{equation}
\beta_1 \beta_2 \beta_3+  \beta_1 \gamma_{2,3}+ \gamma_{1,2} \beta_3 \ .
\end{equation}
This is $Q_3$. So
\begin{eqnarray}
Y_5&=&(\beta_0 \beta_1 \beta_2 \beta_3+ \beta_0 \beta_1 \gamma_{2,3}+ \beta_0 \gamma_{1,2} \beta_3+ \gamma_{0,1} \beta_2\beta_3+\gamma_{0,1}\gamma_{2,3})\ Y_1 \nonumber \\
&+&\gamma_{-1,0}(\beta_1 \beta_2 \beta_3+  \beta_1 \gamma_{2,3}+ \gamma_{1,2} \beta_3)\ Y_0 \ .
\end{eqnarray}

\section{The second solution via the second-order differential equation}

\subsection{The extended cauchy-integral method}

The close connection between differential equations and recurrence relations
enables one to use solutions to differential equations to generate solutions
to the corresponding recurrence relations. Functions of the hypergeometric type
are of particular interest in that they obey a second order differential equation
in the continuous independent variable, and a difference equation in any one
of its parameters.

In the following sections we use the extended Cauchy-integral
method not only to obtain a second solution to the differential equation but to provide
as well a second solution to the difference equation obeyed by one of the parameters.
The results obtained using the Cauchy-integral method are then also shown to follow
from d'Alembert's reduction of order method.\footnote{For a representation of hypergeometric
second-kind solutions using a Rodrigues-type formula, see Area et al.\ \cite{Area}.}
We illustrate these methods by considering the differential equation for the
confluent hypergeometric function:
\begin{equation}\label{dech}
xy'' + (b-x)y' -ay = 0
\end{equation}
in which $a=-N$, a non-positive integer, and $b=n+1$, a positive integer. Although the
analysis that follows requires $b=n+1$, our original interest in this choice
of parameters was the observation that the two standard solutions, $\,_1\!F_1(a;b;x)$
and $U(a,b,x)$, are no longer independent provided only that $a=-N$, in which case
$U(-N,b,x) = (-1)^{N}(b)_N\!\ \!_1\!F\!_1(-N;b;x)$ (see DLMF \cite[Eqs. 13.2.7, 13.2.10 and
13.2.34]{DLMF}).\footnote{We use the Pochhammer symbol defined by $(b)_N=\Gamma(b+N)/\Gamma(b)$.}

Equation (\ref{dech}) has a polynomial solution defined by the confluent hypergeometric function
\begin{equation}\label{Phi}
\Phi(N,n,x) \equiv \,_1\!F_1(-N;n+1;x)  = \sum_{k=0}^{N}\frac{(-N)_{k}}{(n+1)_{k}k!}x^{k},
\end{equation}
which constitutes a first solution of the differential equation in $x$ and a first solution
of the difference equation in the first parameter, $N$:
\begin{eqnarray}
&&(N+n+1)\Phi(N+1,n,x) - (2N+n+1-x)\Phi(N,n,x) +N\Phi(N-1,n,x)\nonumber\\
&&=0 . \label{dephi}
\end{eqnarray}
We derive  a polynomial solution to this equation that is linearly independent of the function
$\Phi(N,n,x)$.

Following Nikiforov and Uvarov \cite[\S 11, p.\,97, Eq. (4)]{NU},
a second linearly independent solution to Eq. (\ref{dech}) is given by
the extended Cauchy integral:
\begin{equation}\label{cipsi}
\Psi(N,n,x) = \frac{1}{\rho(x)}\int_{0}^{\infty}\frac{\rho(s)\Phi(N,n,s)}{s-x}\,ds
\end{equation}
in which the weight function $\rho(x)$ = $e^{-x}x^{b-1}$ is, for the differential equation (\ref{dech}),
a solution of the equation $(x\rho(x))' = (b-x)\rho(x)$.

We next show that $\Psi(N,n,x)$ as defined in (\ref{cipsi}) with $b=n+1$ obeys the same difference equation as $\Phi(N,n,x)$:
\begin{samepage}
\footnotesize{
\begin{eqnarray}\label{difeqpsi}
&&(N+n+1)\Psi(N+1,n,x) -(2N+n+1-x)\Psi(N,n,x) +N\Psi(N-1,n,x)\nonumber\\
&& = x^{-n}e^{x}\int_{0}^{\infty}ds\frac{e^{s}s^{n}}{s-x}\times\nonumber\\
&&\big[(N+n+1)\Phi(N+1,n,s) -(2N+n+1-x)\Phi(N,n,s) +N\Phi(N-1,n,s)\big]\,ds\nonumber\\
&& = x^{-n}e^{x}\int_{0}^{\infty}ds\frac{e^{s}s^{n}}{s-x}\times\nonumber\\
&&\Big[(N+n+1)\Phi(N+1,n,s) -(2N+n+1-s)\Phi(N,n,s) +N\Phi(N-1,n,s)\nonumber\\
&& \hspace{4cm}-(s-x)\Phi(N,n,s)\Big]\nonumber\\
&& =  -x^{-n}e^{x}\int_{0}^{\infty}e^{-s}s^{n}\Phi(N,n,s)\,ds
\end{eqnarray}
}
\end{samepage}
in view of (\ref{dephi}). Substituting (\ref{Phi}) in the last integral in (\ref{difeqpsi}) we have
\begin{eqnarray}
\int_{0}^{\infty}e^{-s}s^{n}\Phi(N,n,s)\,ds &=& \int_{0}^{\infty}e^{-s}s^{n}\sum_{k=0}^{N}\frac{(-N)_{k}s^{k}}{(n+1)_{k}k!}ds\nonumber \\
&=& \sum_{k=0}^{N}\frac{(-N)_{k}}{(n+1)_{k}k!}\int_{0}^{\infty}e^{-s}s^{n+k}\,ds\nonumber \\
&=& \sum_{k=0}^{N}\frac{(-N)_{k}\Gamma(n+k+1)}{(n+1)_{k}k!}\nonumber \\
&=& n!\sum_{k=0}^{N}(-1)^{k}\binom{N}{k}\nonumber  \\
&=& n!\,(1-1)^{N}= 0
\end{eqnarray}
for $N=1,2,\ldots .$ Thus, the function $\Psi(N,n,x)$ satisfies the difference equation (\ref{dephi}), i.e.,
\begin{eqnarray}
&&(N+n+1)\Psi(N+1,n,x) -(2N+n+1-x)\Psi(N,n,x) +N\Psi(N-1,n,x) \nonumber \\
&&= 0 \label{depsi}
\end{eqnarray}

Next we write (\ref{cipsi}) in the form
\begin{eqnarray}
\Psi(N,n,x)&=&x^{-n}e^x\int_0^{\infty}\frac{e^{-s}}{s-x}\left[s^n\Phi(N,n,s)-x^n\Phi(N,n,x)\,\right]\,ds\nonumber\\
&&\ \ +\ \Phi(N,n,x)\int_{-x}^{\infty}\frac{e^{-s}}{s}ds \ . \label{Psi-sol0}
\end{eqnarray}
Inserting the polynomial expression (\ref{Phi}) for $\Phi$, the first of the two integrals in
Eq.\,(\ref{Psi-sol0}) is
\begin{samepage}
\begin{eqnarray}
\overline{\Psi}(N,n,x) &\equiv& \sum_{k=0}^{N}\frac{(-N)_{k}}{(n+1)_{k}k!}\int_{0}^{\infty}\frac{e^{-s}}{(s-x)}\big[s^{n+k} - x^{n+k}\big]ds\nonumber \\
 &=& \sum_{k=0}^{N}\frac{(-N)_{k}}{(n+1)_{k}k!}\sum_{m=0}^{n+k-1}(n+k-1-m)!\,x^{m}
\end{eqnarray}
\end{samepage}
while in the last term of Eq. (\ref{Psi-sol0}) the exponential-integral function
\begin{equation}
\mbox{Ei}\left(x\right)=-\int_{-x}^{\infty }\frac{e^{-s}}{s}\,ds\
\end{equation}
appears. There results\footnote{%
It is worth noting that an analogous result exists for the Legendre polynomials: $\textrm{Q}_{n}(x) =
-W_{n-1}(x) +\textrm{P}_{n}(x)\ln{\sqrt{(1+x)/(1-x)}}$ in which $W_{n-1}$ is a polynomial of order $n-1$, as shown in Erd\'{e}lyi \cite[\S 3.6.2, Eq.(26)]{Erdelyi}, and reflecting the
natural separation of second solutions to homogeneous second-order hypergeometric differential equations into a so-called "polynomial" part and a "logarithmic part".}
\begin{eqnarray}\label{psibar}
\Psi(N,n,x) &=& \overline{\Psi}(N,n,x)\frac{e^{x}}{x^{n}} - \Phi \left( N,n,x\right) \mbox{Ei}(x)\nonumber \\
&=&\frac{n!}{(N+n)!}\,P( N,n,x)\frac{e^{x}}{x^{n}}
-\Phi \left( N,n,x\right) \mbox{Ei}(x) \label{Psi}
\end{eqnarray}
where the polynomial $P\left( N,n,x\right) $ is
\begin{eqnarray}
&&P\left( N,n,x\right) =\frac{(N+n)!}{n!}%
\sum_{k=0}^{N}\sum_{m=0}^{n+k-1}\frac{(-N)_{k}(n+k-1-m)!}{(n+1)_{k}\ k!}x^{m} \
\label{P1}
\end{eqnarray}
and
\begin{equation}
\overline{\Psi}(N,n,x)=\frac{n!}{(N+n)!}\,P( N,n,x) \ .\label{Psi-bar}
\end{equation}
The normalization of the polynomial $P(N,n,x)$ has been chosen to make the coefficient of $x^{N+n-1}$
be $(-1)^N$.  It then turns out that all the coefficients of the powers of $x$ are integers.

Since both $\Phi(N,n,x)$ and $\Psi(N,n,x)$ satisfy the difference equation (\ref{dephi}), it follows that  $\overline{\Psi}(N,n,x)$ also satisfies this equation, as factors independent of  $N$, such as $e^{x}$
and Ei$(x)$ in Eq.\,(\ref{psibar}), do not modify the difference equation. Moreover,
we can show that $\overline{\Psi}(N,n,x)$ is linearly independent of $\Phi(N,n,x)$:
Multiplying (\ref{dephi}) by $\overline{\Psi}(N,n,x)$ and (\ref{depsi})
written for $\overline{\Psi}$ by $\Phi(N,n,x)$ and subtracting, we have
\begin{equation}\label{cass}
(N+n+1)\mathcal{C}(N+1) = N\mathcal{C}(N)
\end{equation}
where $\mathcal{C}$ is the Casoratian:
\begin{equation}
\mathcal{C}(N) = \Phi(N,n,x)\overline{\Psi}(N-1,n,x) - \overline{\Psi}(N,n,x)\Phi(N-1,n,x) \ .
\end{equation}
From (\ref{cass})
\begin{equation}
\mathcal{C}(N+1) = \frac{N}{N+n+1}\mathcal{C}(N)
\end{equation}
from which
\begin{equation}
\mathcal{C}(N+1) = \prod_{k=1}^{N}\Big(\frac{k}{k+n+1}\Big)\mathcal{C}(1) \ .
\end{equation}
From (\ref{Phi}) and (\ref{psibar})
\begin{eqnarray}
\Phi(0,n,x) &=& 1\\
\Phi(1,n,x) &=& 1-\frac{x}{n+1}\\
\overline{\Psi}(0,n,x) &=& \sum_{m=0}^{n-1}(n-1-m)!\,x^{m}\\
\overline{\Psi}(1,n,x) &=& \int_{0}^{\infty}\frac{e^{-s}}{s-x}\bigg[s^{n}\Big(1-\frac{s}{n+1}\Big) - x^{n}\Big(1-\frac{x}{n+1}\Big)\bigg]\nonumber \\
&=& \sum_{m=0}^{n-1}(n-1-m)!\,x^{m} - \frac{1}{n+1}\sum_{m=0}^{n}(n-m)!\,x^{m}\nonumber \\
&=& \sum_{m=0}^{n-1}(n-1-m)!\,x^{m} - \frac{1}{n+1}\sum_{m=-1}^{n-1}(n-1-m)!\,x^{m+1}\nonumber \\
&=& \Big[1-\frac{x}{n+1}\Big]\sum_{m=0}^{n-1}(n-1-m)!\,x^{m}\, -\, \frac{n!}{n+1} \ .
\end{eqnarray}
We then have
\begin{equation}
\mathcal{C}(1) = \Phi(1,n,x)\overline{\Psi}(0,n,x)-\overline{\Psi}(1,n,x)\Phi(0,n,x)=\frac{n!}{n+1}
\end{equation}
from which
\begin{equation}
\mathcal{C}(N+1) = \frac{N!(n!)^{2}}{(N+n+1)!} \ ,
\end{equation}
thus proving that the polynomial $\frac{n!}{(N+n)!}\,P(N,n,x)$ and $\Phi(N,n,x)$ are linearly independent
solutions of the difference equation (\ref{dephi}).

In order to find the coefficients of the powers of $x$ in $P(N,n,x)$, we interchange the order of
summations in (\ref{P1}):

\begin{equation}\label{dubsum}
\sum_{k=0}^{N}\sum_{m=0}^{n+k-1} = \sum_{m=0}^{n-1}\sum_{k=0}^{N}  + \sum_{m=n}^{N+n-1}\sum_{k=m-n+1}^{N} \ .
\end{equation}
From the first double sum on the right-hand side we have,
writing $(n+k-1-m)! = (n-m)_{k}(n-m-1)!$,
\begin{equation}
\sum_{m=0}^{n-1}x^{m}(n-m-1)!\sum_{k=0}^{N}\frac{(-N)_{k}(n-m)_{k}}{(n+1)_{k}k!}
\end{equation}
in which we can use Gauss' formula $\,_2F_1(a,b;c;1) = \frac{\Gamma(c)\Gamma(c-a-b)}{\Gamma(c-a)\Gamma(c-b)}$ to write
\begin{equation}
\sum_{k=0}^{N}\frac{(-N)_{k}(n-m)_{k}}{(n+1)_{k}k!} = \,_2F_1(-N,n-m;n+1;1) = \frac{n!(N+m)!}{(N+n)!m!} \ .
\end{equation}
Thus the first double sum on the right-hand side of (\ref{dubsum}) gives
\begin{equation}
\frac{n!}{(N+n)!}\sum_{m=0}^{n-1}\frac{(N+m)!(n-m-1)!}{m!}x^{m} \ .
\end{equation}
From (\ref{P1}) and (\ref{dubsum}) there results
\begin{samepage}
\begin{eqnarray}
&&P\left( N,n,x\right) \nonumber\\
&&=\sum_{m=0}^{n-1}\left[\frac{\left( N+m\right) !\left( n-m-1\right)
!}{m!}\right]x^{m}\nonumber\\
&&-\,x^{n}\sum_{m=0}^{N-1}%
\left[\sum_{k=0}^{N-m-1}\frac{N!}{\left( N-k-m-1\right) !}%
\frac{(N+n)!}{\left( n+k+m+1\right) !}\frac{(-1)^k k!}{\left( k+m+1\right) !}\right]\left(
-x\right) ^{m}. \label{P}\nonumber\\
\end{eqnarray}
\end{samepage}
In Appendix \ref{resum} we show how to simplify the inner sum in the second term.

We give here a few explicit cases for the polynomial $P(N,n,x):$
\begin{equation}
\def\arraystretch{1.0}
\begin{array}{cc}
n & N=0 \\
0 & 0 \\
1 & 1 \\
2 & \,x\p 1 \\
3 & \,x^2\p \,x\p 2 \\
4 & \,x^3\p \,x^2\p 2 \,x\p 6 \\
5 & \,x^4\p \,x^3\p 2 \,x^2\p 6 \,x\p 24 \\
6 & \,x^5\p \,x^4\p 2 \,x^3\p 6 \,x^2\p 24 \,x\p 120
\end{array}
\end{equation}

\begin{equation}
\def\arraystretch{1.0}
\begin{array}{cc}
n & N=1 \\
0 & \m 1 \\
1 & \m \,x\p 1 \\
2 & \m \,x^2\p 2 \,x\p 1 \\
3 & \m \,x^3\p 3 \,x^2\p 2 \,x\p 2 \\
4 & \m \,x^4\p 4 \,x^3\p 3 \,x^2\p 4 \,x\p 6 \\
5 & \m \,x^5\p 5 \,x^4\p 4 \,x^3\p 6 \,x^2\p 12 \,x\p 24 \\
6 & \m \,x^6\p 6 \,x^5\p 5 \,x^4\p 8 \,x^3\p 18 \,x^2\p 48 \,x\p 120
\end{array}
\end{equation}

\begin{equation}
\def\arraystretch{1.1}
\begin{array}{cc}
n & N=2 \\
0 & \,x\m 3 \\
1 & \,x^2\m 5 \,x\p 2 \\
2 & \,x^3\m 7 \,x^2\p 6 \,x\p 2 \\
3 & \,x^4\m 9 \,x^3\p 12 \,x^2\p 6 \,x\p 4 \\
4 & \,x^5\m 11 \,x^4\p 20 \,x^3\p 12 \,x^2\p 12 \,x\p 12 \\
5 & \,x^6\m 13 \,x^5\p 30 \,x^4\p 20 \,x^3\p 24 \,x^2\p 36 \,x\p 48 \\
6 & \,x^7\m 15 \,x^6\p 42 \,x^5\p 30 \,x^4\p 40 \,x^3\p 72 \,x^2\p 144 \,x\p 240
\end{array}
\end{equation}

\begin{equation}
\def\arraystretch{1.0}
\begin{array}{cc}
n & N=3 \\
0 & \m \,x^2\p 8 \,x\m 11 \\
1 & \m \,x^3\p 11 \,x^2\m 26 \,x\p 6 \\
2 & \m \,x^4\p 14 \,x^3\m 47 \,x^2\p 24 \,x\p 6 \\
3 & \m \,x^5\p 17 \,x^4\m 74 \,x^3\p 60 \,x^2\p 24 \,x\p 12 \\
4 & \m \,x^6\p 20 \,x^5\m 107 \,x^4\p 120 \,x^3\p 60 \,x^2\p 48 \,x\p 36 \\
5 & \m \,x^7\p 23 \,x^6\m 146 \,x^5\p 210 \,x^4\p 120 \,x^3\p 120 \,x^2\p 144 \,x\p 144 \\
6 & \m \,x^8\p 26 \,x^7\m 191 \,x^6\p 336 \,x^5\p 210 \,x^4\p 240 \,x^3\p 360 \,x^2\p 576 \,x\p 720
\end{array}
\end{equation}

\begin{equation}
\def\arraystretch{1.0}
\begin{array}{cc}
n & N=4 \\
0 & \,x^3\m 15 \,x^2\p 58 \,x\m 50 \\
1 & \,x^4\m 19 \,x^3\p 102 \,x^2\m 154 \,x\p 24 \\
2 & \,x^5\m 23 \,x^4\p 158 \,x^3\m 342 \,x^2\p 120 \,x\p 24 \\
3 & \,x^6\m 27 \,x^5\p 226 \,x^4\m 638 \,x^3\p 360 \,x^2\p 120 \,x\p 48 \\
4 &\,x^7\m 31 \,x^6\p 306 \,x^5\m 1066 \,x^4\p 840 \,x^3\p 360 \,x^2\p 240 \,x\p 144 \\
5 & \,x^8\m 35 \,x^7\p 398 \,x^6\m 1650 \,x^5\p 1680 \,x^4\p 840 \,x^3\p 720 \,x^2\p 720 \,x\p 576 \\
6 & x^9\m 39 x^8\p 502 x^7\m 2414 x^6\p 3024 x^5\p 1680 x^4\p 1680 x^3\p 2160 x^2\p 2880 x\p 2880
\end{array}
\end{equation}

The coefficients of $x$ in the polynomial $P\left( N,n,x\right) $ up to the
power $x^{n-1}$\ are all positive and contain relatively simple (factorial)
factors, while those for powers $x^{n}$ up to the highest power $x^{N+n-1}$
have oscillating signs and some may have very high prime number
factors, much larger than $N+n-1$, so that they will not reduce to simple
factorials. Rather, the coefficients for powers at and above $x^n$ involve
harmonic sums.\footnote{%
These properties of the coefficients apply even more generally to the second
solution polynomials allied with the full hypergeometric functions $%
\,_2F_1(-N,b;c;x)$, but we leave the explicit derivation
of these polynomials to the especially engaged reader.} (See Appendix \ref{resum}.)

The function $\Phi(N,n,x)$ considered here is, apart from a normalization factor, the
well-known associated Laguerre polynomial $L_N^{(n)}(x)$
(DLMF \cite[Eq.\,18.5.12]{DLMF}). We therefore define a suitably normalized associated
Laguerre polynomial of the first kind with
\begin{equation}
\overline{\strut L}_N^{(n)}(x)\equiv\frac{N!\,n!}{(N+n)!}L_N^{(n)}(x)=\Phi(N,n,x)
\end{equation}
and an associated Laguerre function of the second kind with
\begin{equation}
\overline{\overline{\strut L}}_{N}^{(n)}(x)\equiv \frac{n!}{(N+n)!}\,\frac{1}{x^n}P(N,n,x) \ .
\end{equation}

Both of these Laguerre functions, $\overline{\strut L}_N^{(n)}(x)$ and  $\overline{\overline{\strut L}}_{N}^{(n)}(x)$,  satisfy recurrence relations
as given by DLMF \cite[Eqs. 13.3.1,13.3.2]{DLMF}:
\begin{equation}
(b-a)\,_1\!F_1(a-1;b;x)+(2a-b+x)\,_1\!F_1(a;b;x)-a\,_1\!F_1(a+1;b;x)=0%
\label{recurs-a} 
\end{equation}
\begin{equation}
b(b-1)\,_1\!F_1(a;b-1;x)+b(1-b-x)\,_1\!F_1(a;b;x)+x(b-a)\,_1\!F_1(a;b+1;x)=0 \label{recurs-b}
\end{equation}
where $\,_1\!F_1(a;b;x)$ is a confluent hypergeometric function with $a=-N$, $b=n+1$.

We have, for either $\overline{\strut L}_N^{(n)}(x)$ or for $\overline{\overline{\strut L}}_N^{(n)}(x),$
\begin{equation}
(n+1+N)\,\overline{\overline{\strut L}}_{N+1}^{(n)}-(2N+n+1-x)\,\overline{\overline{\strut L}}_{N}^{(n)}+N\, \overline{\overline{\strut L}}_{N-1}^{(n)}=0 \label{recurs-N}
\end{equation}
as well as
\begin{equation}
n(n+1)\,\overline{\overline{\strut L}}_N^{(n-1)}-(n+1)(n+x)\,\overline{\overline{\strut L}}_N^{(n)}+x(n+1+N)\,\overline{\overline{\strut L}}_N^{(n+1)}=0 \ .
\end{equation}

\subsection{D'Alembert's reduction of order method}

We now use d'Alembert's reduction-of-order method to generate the second solution
to the confluent hypergeometric recurrence relation found in the previous section
using the extended Cauchy-integral method.
We have, for any second-order
homogeneous linear differential equation, taken in the form
\begin{equation}\label{dif-gen}
y^{\prime \prime }+p(x) y^{\prime }+q(x) y=0
\end{equation}
with a known solution
\begin{equation}
y_{1}=\Phi ( x)\, ,
\end{equation}
a second solution that can be found with the ansatz
\begin{equation}
y_{2}=f(x) \Phi (x) \ .
\end{equation}
Substituting, one finds
\begin{equation}
y_{2}=\Phi \int \frac{1}{\Phi ^{2}}\exp \left( -{\textstyle\int p\, dx} \right) dx \ .
\label{inv-phi2}
\end{equation}

The Wronskian $W\left( y_{1},y_{2}\right)\equiv y_1{y^{\prime}}_2-y_2{%
y^{\prime}}_1 =\exp \left( -\int p\,dx\right) $, which implies that in the
range of $x$ for which $\int p\,dx$ is not infinite, the two solutions are
independent.

The integral in Eq.\,(\ref{inv-phi2}) looks difficult in cases in which $%
\Phi(x)$ is not simple. However, this indefinite integral marvelously
simplifies\footnote{%
Although many such intriguing integrals can be generated, such as the
Legendre case
$Q_n(x)=P_n(x){\begingroup\textstyle\int\endgroup}^x dx'\,(1-{x'}^2)^{-1}/[P_n(x')]^{\,2},$ 
only a few non-trivial examples appear in the commonly-used reference compilations.
Gradshteyn and Ryzhik \cite[Eq.\,6.539.1]{Grad}, have
${\begingroup\textstyle\int\endgroup}_a^b dx\,x^{-1}/\left[J_{\nu}(x)\right]^2=(\pi/2)\left[N_{\nu}(b)/J_{\nu}(b)-N_{\nu}(a)/J_{\nu}(a)\right]$, an expression derived by
E.\,von\,Lommel in 1871 and reproduced by Watson \cite[\S\,5.11(3)]{Watson}.
Such relations are easily generated by integrating the identity
$W(y_1,y_2)/y_1^2\equiv (y_2/y_1)^{\prime}$, where the Wronskian of the two independent
solutions $y_1$ and $y_2$ is proportional to $\exp{(-{\begingroup\textstyle\int\endgroup} p\,dx)}$.}
when $\Phi$ is a solution to the second-order equation (\ref
{dif-gen}).
In the case of rational integrands, perhaps with transcendental  arguments and algebraic factors, algorithms for doing
such indefinite integrals now exist.\footnote{For a description of these
methods, see Bronstein \cite{Bronstein} and also Geddes et al. \cite{Geddes}. Many have been implemented in
a variety of symbolic manipulation programs, notably {\it Mathematica}
and {\it Maple}%
.}
Note, however, that the simplification that occurs in the
integration specified in Eq.\,(\ref{inv-phi2}) is delicate, in that the integers that appear
in the polynomials in Eq.\ (\ref{inv-phi2}) must be precisely those in the polynomial $\Phi(x)$. Slight deviations can cause
an explosion of extra terms in the resultant integral.

If we apply Eq.\,(\ref{inv-phi2}) to our confluent hypergeometric
differential equation (\ref{dech}) with $a=-N$ and $b=n+1$,
then we can take
$y_1=\Phi(N,n,x)$
as a polynomial first solution, and have
\begin{equation}\label{px}
p\left( x\right) = \left( n+1\right) /x - 1,
\end{equation}
so that an independent second
solution will be
\begin{equation}
y_{2}(x)=\Phi(x)\int_{-\infty}^x \frac{e^{s}}{\left[\Phi(s)\right]^2}\frac{ds}{s^{n+1}} \ .  \label{inv}
\end{equation}

Performing the integrations (described in Appendix \ref{inv-poly})
we find that this second solution matches that found earlier (Eq.\,(\ref{psibar})):
\begin{eqnarray}
&&y_{2}( -N,n+1,x) \nonumber \\
&=&-\frac{(N+n)!}{N!(n!)^2}\left(\frac{n!}{(N+n)!}P( N,n,x)\frac{e^{x}}{x^{n}}
-\Phi (N,n,x)\mbox{Ei}(x)\right) \ , \label{diff-sol}
\end{eqnarray}
apart from the overall factor.
One can verify this overall factor in the
case of positive integer $n$ by using
\begin{equation}
\lim_{x \to 0}x^n e^{-x}\Phi(x)\int_{-\infty}^x \frac{e^{s}}{\left[\Phi(s)\right]^2}\frac{ds}{s^{n+1}}=
\int_{-\infty}^1\frac{dt}{t^{n+1}}
=-\frac{1}{n}
\end{equation}
and
\begin{equation}
\lim_{x \to 0}x^n e^{-x}\Phi (N,n,x)\mbox{Ei}(x)=\lim_{x \to 0}x^n\mbox{Ei}(x)=0 \ ,
\end{equation}
to find from Eq.\,(\ref{diff-sol}) that
\begin{equation}
P(N,n,0)=N!(n-1)! \ .
\end{equation}
This agrees with the normalization of the $P(N,n,x)$ we selected that has
the coefficient of the highest power of $P(N,n,x)$ to be $(-1)^N$.

\section{Conclusion}

As expected, systematic methods can be developed and applied for finding second solutions to
linear second-order difference equations, analogous to those for differential equations.
We have applied these methods to find a general solution to the confluent hypergeometric recurrence relations Eqs.\,(\ref{recurs-a}, \ref{recurs-b}) in
the degenerate case $a=-N,\ b=n+1$ ($N\ge 0$ and $n\ge 0$).
The second solution to these recurrence relations is proportional to
$(n!/(N+n)!)x^{-n}P(N,n,x)$, where the polynomial $P(N,n,x)$ is given by Eq.\,(\ref{P}).
In particular, the second solution to just the recurrence relation of
Eq.\,(\ref{recurs-N}) is proportional to the polynomial
${\overline\Psi(N,n,x)}=(n!/(N+n)!)P(N,n,x).$

Curiously, the closed-form second solution to the confluent hypergeometric
differential equation in the degenerate case when the first parameter $a$ in $%
\,_1\!F_1(a;b;x)$ takes the value of a non-positive integer
and the second parameter $b$ is an integer greater than zero is not yet found
in standard references. For example, the DLMF gives instead an infinite Laurent power-series
representation (see DLMF \cite[Eq.\,13.2.31]{DLMF}), a result which we
reconstruct in Appendix \ref{ser-sol}.  In Appendix \ref{compare},
we show that the DLMF expression matches our closed-form solution Eq.~\eqref{psibar}.

\section*{Acknowledgements}
Both authors gratefully acknowledge the support of The George Washington University through its Physics Department, and the second author thanks the Arizona State University Physics Department for the effortless accessibility of online materials.



\newpage

\appendix

\section{Re-summing in the confluent hypergeometric polynomial of the second kind}
\label{resum}

Our confluent hypergeometric second solution polynomial is given by Eq.\,\ref{P} as:
{\small
\begin{eqnarray}
&&P\left( N,n,x\right)  \nonumber \\
&&=\sum_{m=0}^{n-1}\frac{\left( N+m\right) !\left(
n-m-1\right) !}{m!}x^{m}-x^n\times\nonumber \\
&&\sum_{m=0}^{N-1}\sum_{k=0}^{N-m-1}\left( -1\right) ^{k}\frac{N!}{%
\left( N-k-m-1\right) !}\frac{\left( N+n\right) !}{\left( n+k+m+1\right) !}%
\frac{k!}{\left( k+m+1\right) !}\left( -x\right) ^{m}\ .\ \ \ \ \ \label{P-apen}
\end{eqnarray}
}
The expression for the coefficients in the inner sum of the second term, which we write as
{\footnotesize
\begin{eqnarray}
&&c\left( N,n,m\right)  \nonumber \\
&&=(-1)^{m+1}\sum_{k=0}^{N-m-1}\left( -1\right) ^{k}\frac{N!}{%
\left( N-k-m-1\right) !}\frac{\left( N+n\right) !}{\left( n+k+m+1\right) !}%
\frac{k!}{\left( k+m+1\right) !}
 \nonumber\\
&&=\left( -1\right) ^{m+1}\frac{N!}{\left( N-m-1\right) !}\frac{\left(
n+N\right) !}{\left( m+n+1\right) !\left( m+1\right) !}
\,_{3}F_{2}\left(
1,1,-N+m+1;2+m,2+m+n;1\right)  \nonumber
\end{eqnarray}
}
\vspace{-18 mm}
\begin{equation}
\label{c-apend}
\end{equation}
can be simplified considerably.

First, we re-express the $\,_3F_2$ hypergeometric polynomial in terms of an
integral over an $\,_2F_1$ hypergeometric function (DLMF \cite[16.5.2]{DLMF})
\begin{eqnarray}
&&\,_{3}F_{2}\left( a_{1},a_{2},c;b,d;z\right)  \nonumber \\
&=&\frac{\Gamma \left( d\right) }{\Gamma \left( c\right) \Gamma \left(
d-c\right) }\int_{0}^{1}t^{c-1}\left( 1-t\right) ^{d-c-1}\,_{2}F_{1}\left(
a_{1},a_{2},b,zt\right) dt \ .
\end{eqnarray}
In our case,
\begin{eqnarray}
&&\,_{3}F_{2}\left( -N+m+1,1,1;m+n+2,m+2;1\right)  \nonumber\\
&=&\left( m+1\right) \int_{0}^{1}\left( 1-t\right) ^{m}\,_{2}F_{1}\left(
-N+m+1,1;m+n+2;t\right) dt \ .
\end{eqnarray}
In turn, the $\,_2F_1$ can be written as an integral (DLMF \cite[15.6.1]{DLMF}), giving
\begin{samepage}
\begin{eqnarray}
&&\,_{3}F_{2}\left( -N+m+1,1,1;m+n+2,m+2;1\right) \nonumber \\
&=&\left( m+1\right) \left( m+n+1\right) \int_{0}^{1}\left( 1-t\right)
^{m}\left( \int_{0}^{1}\left( 1-s\right) ^{m+n}\left( 1-st\right)
^{N-m-1}ds\right) \,dt \ . \nonumber \\
\end{eqnarray}
\end{samepage}
Now the trick for re-summing in our coefficients is to expand the integrand
factor $\left( 1-st\right)^{N-m-1} $ not in $s$, but in $1-s\equiv u$:
\begin{samepage}
{\footnotesize
\begin{eqnarray}
&&\,_{3}F_{2}\left( -N+m+1,1,1;m+n+2,m+2;1\right) \nonumber \\
&=&\left( m+1\right) \left( m+n+1\right) \int_{0}^{1}\left( 1-t\right)
^{m}\left( \int_{0}^{1}u^{m+n}\left( 1-t+tu\right) ^{N-m-1}du\right) \,dt \nonumber  \\
&=&\left( m+1\right) \left( m+n+1\right) \sum_{k=0}^{N-m-1}
\left(
\begin{array}{c}
N-m-1\\ k
\end{array}
\right)
\int_{0}^{1}\left( 1-t\right) ^{N-1-k}t^{k}\left(
\int_{0}^{1}u^{m+n}u^{k}du\right) \,dt \nonumber \\
&=&\left( m+1\right) \left( m+n+1\right) \sum_{k=0}^{N-m-1}
\left(
\begin{array}{c}
N-m-1\\ k
\end{array}
\right)
\frac{1}{m+n+k+1}\int_{0}^{1}t^{k}\left( 1-t\right) ^{N-1-k}\,dt \nonumber \\
&=&\left( m+1\right) \left( m+n+1\right) \sum_{k=0}^{N-m-1}\frac{\left(
N-m-1\right) !}{\left( N-m-1-k\right) !k!}\frac{k!\left( N-k-1\right) !}{N!}%
\frac{1}{m+n+k+1}\nonumber  \\
&=&\frac{\left( m+1\right) \left( m+n+1\right) }{N!}\left( N-m-1\right)!
\sum_{k=0}^{N-m-1}\frac{\left( N-k-1\right) !}{\left( N-m-1-k\right) !}%
\frac{1}{m+n+k+1}\nonumber
\end{eqnarray}
}
\end{samepage}
where we used
\begin{equation}
\int_{0}^{1}t^{a}\left( 1-t\right) ^{b}dt=\frac{\Gamma \left( a+1\right)
\Gamma \left( b+1\right) }{\Gamma \left( a+b+2\right) }.
\end{equation}

We have arrived at an alternate and simpler expression for our coefficients:
\begin{eqnarray}
&&c\left( N,n,m\right) \nonumber  \\
&=&\left( -1\right) ^{m+1}\frac{\left( n+N\right) !}{\left( n+m\right) !m!}%
\sum_{k=0}^{N-m-1}\frac{\left( N-1-k\right) !}{\left( N-1-k-m\right) !}\frac{%
1}{n+m+k+1}\nonumber\\
&=&\left( -1\right) ^{m+1}\frac{\left( n+N\right) !}{%
\left( n+m\right) !}\sum_{k=n+m+1}^{N+n}
\left(
\begin{array}{c}
N+n+m-k\\ m
\end{array}
\right)
\frac{1}{k}\nonumber \\
&=&(-1) ^{m+1}\frac{( n+N)!}{(n+m)!m!}%
\sum_{k=n+m+1}^{N+n}(N+n+1-k)_m\frac{1}{k} \ . \label{c-harm}
\end{eqnarray}

As a check, note that if $m=N-1$ (corresponding to the highest power of $x$ in our polynomial $P(N,m,x)$),
\begin{equation*}
c\left( N,n,N-1\right) =\left( -1\right) ^{N} \ .
\end{equation*}
The case $m=0$ is particularly simple, and shows that harmonic numbers enter
these coefficients.
\begin{equation}
c\left( N,n,0\right) =-\frac{\left( n+N\right) !}{n!}\sum_{k=n+1}^{N+n}%
\frac{1}{k} \ .
\end{equation}
The case $m=1$ leads to
\begin{equation}
c\left( N,n,1\right)
=\frac{\left( n+N\right) !}{\left( n+1\right) !}\left( \left( N+n+1\right)
\sum_{k=n+2}^{N+n}\frac{1}{k}-\left( N-1\right) \right) \ .
\end{equation}
The Pochhammer factor $(N+n+1-k)_m$ in Eq.\,(\ref{c-harm}) is a
polynomial of degree $m$ in the
summation variable $k$. The contribution to $c(N,n,m)$ from this polynomial expanded in
powers of $k$ will be a harmonic sum
from the $k^0$ term while the terms for powers of  $k$ from $1$ to $m$ will, after
canceling the denominator $k$, lead to a polynomial in $N$ and $n$.

\section{Handling integrals over inverse polynomials}
\label{inv-poly}

The algorithms implemented in presently-available symbolic programs {such as \it Mathematica} and {\it Maple}
can solve a variety of indefinite integrals over integrands containing ratios of polynomials, multiplied by
algebraic and transcendental functions.  These methods do not require knowledge of the roots of the
denominator polynomial.\footnote{Finding exact roots of arbitrary polynomials of degree higher than four
would have been prohibitive.} In our case, Eq.\,(\ref{inv}), we have the inverse of a polynomial squared, namely $[\Phi]^2$,
together with a weight factor. Typically, the first step
is to simplify the integrand by removing the double poles arising from the zeros of the polynomial in the denominator, making the polynomial `square-free'. This can be done as follows. First, note that
\begin{equation}
\left(\frac{f}{\Phi}\right)^\prime = \frac{f'}{\Phi}-\frac{f}{\Phi^2}\Phi^\prime\ .\label{Phi2}
\end{equation}
Here, primes denote a derivative. Now we use B{\'e}zout's identity
(\cite[A.2, p.231]{Dave})
that for any two polynomials $a$ and $b$, there exist polynomials $s$ and $t$ such that
\begin{equation}
sa+tb=gcd(a,b)
\end{equation}
where $gcd(a,b)$ is the greatest common divisor of $a$ and $b$, and the degree of the polynomial $s$ is less than
$deg(b)-deg(gcd(a,b))$ and the degree of $t$ is less than $deg(a)-deg(gcd(a,b))$. There are simple methods
going back to the Babylonians \cite{Knuth}
that extract the greatest common divisors of a pair of integers, requiring only a
sequence of subtractions. These methods extend to polynomials.
Because $\Phi$ and $\Phi^\prime$ are relatively
prime, we have
\begin{equation}
s\Phi+t\Phi^\prime=1\label{Syl}
\end{equation}
where the degree of $s$ is $N-2$, while the degree of $t$ is $N-1$. As an identity in the independent variable of these polynomials, there will be $2N-1$ relations to solve for the $2N-1$ coefficients in the polynomials $s$ and $t$. (Equation (\ref{Syl}) is often expressed
by a $(2N-1)\times(2N-1)$ `Sylvester' matrix times a column vector formed from
the coefficients in the polynomial $s$ and $t$, equal to a column vector $(1,0,0,\cdots)^T$.)  The relations are solvable for a given $N$ by Gaussian reduction; with some effort, it is also possible to construct the polynomials
$s$ and $t$ for arbitrary $N$.\footnote{\label{footst}The polynomials $s$ and $t$, solutions to Eq.\,(\ref{Syl}), are given by
\begin{eqnarray*}
&&s\left( N,n,x\right) =1-\frac{\left( N+n\right) !}{\left( N-1\right) !\left(
n+1\right) !}\\
&&+\frac{\left( N+n\right) !}{\left( N-1\right) !n!}%
\sum_{p=1}^{N-2}\frac{(-N+1)_p}{(n+2)_p}\big( 1-%
\,_3F_2\left( \left[ -N+p+1,p,1\right] ,\left[ n+2+p,-N+1\right]
,1\right) \big)\frac{x^{p}}{p\cdot p!}
\end{eqnarray*}
\begin{eqnarray*}
&&t(N,n,x)=-\frac{\left( N+n\right) !}{N!n!}\sum_{p=0}^{N-1}\frac{(-N+1)_p}{(n+2)_p}\,_3F_2\left( \left[ -N+p+1,p,1\right] ,%
\left[ n+2+p,-N+1\right] ,1\right) \frac{x^{p}}{p!}
\end{eqnarray*}
}

With Eq.\,(\ref{Syl}), $\Phi^{^\prime}$ can be replaced in Eq.\,(\ref{Phi2}) to obtain
\begin{equation}
\left(\frac{f}{\Phi}\right)^\prime =\frac{f^\prime}{\Phi}-\frac{f}{t\Phi^2}+\frac{sf}{t\Phi} \ .
\end{equation}
If we take $f=tg$ and integrate, we will have the identity
\begin{equation}
\int{\frac{g}{\Phi^2}}dx=\int{\frac{(tg)^\prime+sg}{\Phi}}dx-\frac{tg}{\Phi} \ .
\end{equation}
Now for our integral, $g=\exp{(-\int{p\,dx})}$, which satisfies $g'=-pg$, we have
\begin{equation}
\Phi\int{\frac{g}{\Phi^2}}dx=\Phi\int{\frac{t^\prime-pt+s}{\Phi}}gdx-tg \label{one-phi}
\end{equation}
in which $p=p(x)$ is given in (\ref{px}).

The next step in the commonly-used symbolic programs is to first ensure, by repeated subtraction, that
the numerator polynomial, which we call $A$, is of degree less than that of the square-free denominator polynomial, which
we call $B$, and that $A$ and $B$ have no common polynomial divisors. Then the programs apply the Rothstein-Trager method
 and its improvements \cite{Bronstein,Geddes} to
express the integral as a particular sum over the roots $\{z_i\}$ of the resultant constructed from $A$ and $(A-zB^\prime)$.

However, in our case, by employing the solutions given in footnote (\ref{footst}), one can
deduce that
\begin{equation}
t^\prime-pt+s=c\Phi/x
\end{equation}
where $c$ is $( n+1)\left((n+N)!/(n!N!)\right)$.  The denominator $\Phi$ in the integrand of Eq.\,(\ref{one-phi}) is canceled, leaving an exponential integral. One recognizes that this cancellation requires that
the exact solution to $x\Phi^{\prime\prime}+(n+1-x)\Phi^{\prime}+N\Phi=0$ appears in the integrand
denominator. Any other polynomial, even if only slightly different from $\Phi$, will integrate to
a sum over all the zeros of the resultant described above.

\section{The series solution for the confluent hypergeometric function in the degenerate case }
\label{ser-sol}

Because the derivation of the series representation of the second solution for the confluent hypergeometric function
in the degenerate case is not easily found in standard references, we re-derive it here.
Consider the Cauchy integral that produces the standard, regular, first solution to the
confluent hypergeometric equation $zy^{\prime \prime }+\left( b-z\right)
y^{\prime }-ay=0$ when the poles of $\Gamma \left( -s\right) $ are
surrounded: 
\begin{eqnarray}
I\left( z\right) &=&\frac{1}{2\pi i}\oint_{C}\frac{\Gamma \left( a+s\right)
\Gamma \left( -s\right) }{\Gamma \left( b+s\right) }\left( -z\right) ^{s}ds \nonumber \\
&=&\frac{\Gamma \left( a\right) }{\Gamma \left( b\right) }\sum_{k=0}^{\infty }%
\frac{\left( a\right) _{k}}{\left( b\right) _{k}k!}z^{k}\nonumber \\
&=&\frac{\Gamma \left( a\right) }{\Gamma \left( b\right)}\,_1F_1(a;b;z) \ .
\end{eqnarray}
(The Cauchy representation of the hypergeometric function was extensively studied in the
early 1900's by Barnes
\cite{Barnes}.)
We can verify that $I\left( z\right) $ is a solution with more general contours by
\begin{eqnarray}
&&\left( z\frac{d^{2}}{dz^{2}}+\left( b-z\right) \frac{d}{dz}-a\right)
I\left( z\right)  \nonumber \\
&=&\oint_{C'}\frac{\Gamma \left( a+s\right) \Gamma \left( -s+2\right) 
}{\Gamma \left( b+s\right) }\left( -z\right) ^{s-1}ds
+\oint_{C'}b\frac{\Gamma \left(a+s\right) \Gamma \left( -s+1\right) }{\Gamma \left( b+s\right) }\left(
-z\right) ^{s-1}ds \nonumber \\
&&-\oint_{C'}\frac{\Gamma \left( a
+s\right) \Gamma \left( -s+1\right) }{\Gamma \left( b+s\right) }\left( -z\right) ^{s}ds
-\oint_{C'}a\frac{\Gamma \left(
a+s\right) \Gamma \left( -s\right) }{\Gamma \left( b+s\right) }\left(
-z\right) ^{s} ds \nonumber \\
&=&\oint_{C'} \frac{\Gamma \left( a+s\right) }{\Gamma \left( b+s\right) 
}\left(\Gamma \left( -s+1\right) \left( \left( s-1\right) +b\right) \left(
-z\right) ^{s-1}-%
\Gamma \left( -s\right) \left( s+a\right) \left( -z\right) ^{s}\right) ds \nonumber \\
&=&\oint_{C'}\left( \frac{\Gamma \left( a+s\right) }{\Gamma \left(
b+s-1\right) }\Gamma \left( -s+1\right) \left( -z\right) ^{s-1}-\frac{\Gamma
\left( a+s+1\right) }{\Gamma \left( b+s\right) }\Gamma \left( -s\right)
\left( -z\right) ^{s}\right) ds \nonumber\\
&=&\oint_{C'}\left( \frac{\Gamma \left( a+s+1\right) }{\Gamma \left(
b+s\right) }\Gamma \left( -s\right) \left( -z\right) ^{s}-\frac{\Gamma
\left( a+s+1\right) }{\Gamma \left( b+s\right) }\Gamma \left( -s\right)
\left( -z\right) ^{s}\right) ds=0 \ . \nonumber \\
\end{eqnarray}
In the last line, the contour in the first integral was shifted to the right
by one, which will have little affect on the integral if the contour
is over large $s$ where the integrand is negligible
(See Statler \cite[\S 1.8.1]{Lucy}). 

Now consider the integral solution with the contour surrounding all the
poles of the integrand within a large circle. We are particularly interested in the degenerate
case $a=-N,b=n+1$ where $N$ and $n$ are nonnegative integers. Let
\begin{equation*}
I_{C'}(z)=\frac{1}{2\pi i}\oint_{C^{^{\prime }}}\frac{\Gamma \left( -N+s\right) \Gamma
\left( -s\right) }{\Gamma \left( n+1+s\right) }\left( -z\right) ^{s}ds \ .
\end{equation*}
The poles of the integrand in the complex $s$ plane come from the numerator
gamma factors. In the cases for which $s=-n-1,-n-2,\cdots $, the poles of
the numerator are canceled by those in the denominator. There are three
remaining cases for pole contributions: (1) The poles in the integrand
from $s=-n$ to $-1$, which are of order one;
(2) For $s$ from zero to to $N$, the poles of the integrand are of second order due to the
overlap of the poles of the two numerator gamma factors; (3) For $s$ above $%
N,$ the poles in the integrand are order one.

We will use the Cauchy residue theorem in the form 
\begin{equation}
\frac{1}{2\pi i}\oint g\left( z\right) dz=\sum_{k}\lim_{z\rightarrow z_{k}}%
\frac{d^{p_{k}}}{dz^{p_{k}}}\left[ \left( z-z_{k}\right) ^{p_{k}+1}g\left(
z\right) \right] \ ,
\end{equation}
where $g\left( z\right) $ is meromorphic within the contour, and has poles
of order $p_{k}+1$ when $z$ approaches $z_{k}.$

In our case, 
\begin{eqnarray}
g\left( s\right)  &=&\frac{\Gamma \left( a+s\right) \Gamma \left( -s\right) 
}{\Gamma \left( b+s\right) }\left( -z\right) ^{s} \\
\mbox{with }&&  a =-N, \ \ \ \ b =n+1 \ . \nonumber
\end{eqnarray}
Thus
\begin{eqnarray}
&&\frac{1}{2\pi i}\oint g\left( s\right) ds = \nonumber \\
&&\sum_{k=1}^{n}\lim_{s\rightarrow -k}\left[ \left( s+k\right) g\left(
s\right) \right] +\sum_{k=0}^{N}\lim_{s\rightarrow k}\frac{d}{ds}\left(
\left( s-k\right) ^{2}g\left( s\right) \right)\nonumber\\
&&\ \ \ \ \  +\sum_{k=N+1}^{\infty
}\lim_{s\rightarrow k}\left[ \left( s-k\right) g\left( s\right) \right]  \ .
\end{eqnarray}
To calculate the residues, we will be using the identities 
\begin{eqnarray}
\Gamma \left( z\right) \Gamma \left( 1-z\right)  &=&\frac{\pi }{\sin \left(
\pi z\right) } \\
\Gamma \left( 1+z\right)  &=&z\Gamma \left( z\right)  \\
\psi \left( z\right)  &=&\frac{1}{\Gamma \left(
z\right) }\frac{d}{ds}\Gamma \left( z\right)  \\
\psi \left( z\right) -\psi \left( 1-z\right)  &=&-\pi \cot \left( \pi
z\right)  \\
\psi \left( 1+z\right)  &=&\psi \left( z\right) +\frac{1}{z} \\
\psi \left( 1+M\right)  &=&\psi \left( 1\right) +\sum_{k=1}^{M}\frac{1}{k} \\
\psi\left( 1\right)  &=&-\gamma \ .
\end{eqnarray}
There follows 
\begin{eqnarray}
\psi\left( 1+n+s\right)  &=&\psi\left( 1+s\right) +\sum_{l=1}^{n}%
\frac{1}{l}\\
 \mbox{and}&& \nonumber\\
\psi\left( 1+N-s\right)  &=&\psi\left( 1-s\right) +\sum_{l=1}^{N}%
\frac{1}{l} \ .
\end{eqnarray}
We will also employ
\begin{equation}
\cot \left( \pi \left( N-s\right) \right) =-\cot \left( \pi s\right)  \ .
\end{equation}

The single-pole residues for $s=-k+\epsilon $ for small $\epsilon $ and $%
k=1,2,\cdots ,n$ come from those in $\Gamma \left( a+s\right) $ as  
\begin{eqnarray}
&&\lim_{s\rightarrow -k}\left( \left( s+k\right) g\left( s\right) \right)  \nonumber\\
&&=\lim_{s\rightarrow -k}\left( \left( s+k\right) \frac{\Gamma \left(
a+s\right) \Gamma \left( -s\right) }{\Gamma \left( b+s\right) }\left(
-z\right) ^{s}\right)  \nonumber \\
&&=\lim_{\epsilon \rightarrow 0}\epsilon \frac{\Gamma \left( -N-k+\epsilon
\right) \Gamma \left( k-\epsilon \right) }{\Gamma \left( 1+n-k\right) }%
\left( -z\right) ^{-k} \nonumber\\
&&=\left( -1\right) ^{N}\frac{\left( k-1\right) !}{\left( N+k\right) !\left(
n-k\right) !}z^{-k} \ .
\end{eqnarray}
For the double-pole residues, we have 
\begin{eqnarray}
&&\frac{d}{ds}\left( \left( s-k\right) ^{2}g\left( s\right) \right) \nonumber\\
&&=\frac{\left( -1\right) ^{s}z^{s}\Gamma \left( -s\right) \Gamma \left(
-N+s\right) }{\Gamma \left( n+1+s\right) }\nonumber \\
&&\ \ \times\left( s-k\right) ^{2}\left( +%
\psi\left( -N+s\right) -\psi\left( -s\right) -\psi\left(
n+1+s\right) +\ln \left( -z\right) +2\frac{1}{s-k}\right) \nonumber\\
&&=\frac{\left( -z\right) ^{s}\Gamma \left( -s\right) \Gamma \left(
-N+s\right) }{\Gamma \left( 1+n+s\right) }\left( s-k\right) ^{2}\nonumber\\
&&\ \ \times\Big(\psi\left( 1+N-s\right) -\psi\left( 1+s\right) -\psi\left(
1+n+s\right) +\ln \left( -z\right) -2\pi \cot \left( \pi s\right)\nonumber\\
&&\ \ \ \ \ \ \ \ \ \ +2\frac{1}{s-k}\Big) \ 
\end{eqnarray}
in which we have used 
\begin{eqnarray}
\psi\left( -s\right) &=&\psi\left( 1+s\right) +\pi \cot \left(
\pi s\right) \\
\mbox{and}&&\nonumber\\
\psi\left( -N+s\right) &=&\psi\left( 1+N-s\right) -\pi \cot
\left( \pi s\right) \ .
\end{eqnarray}
Near the poles, $s=k+\epsilon ,$ $\left| \epsilon \right| <<1,$ i.e. 
\begin{eqnarray}
-2\pi \cot \left( \pi s\right) +2\frac{1}{s-k} &=&-2\pi \left( \cot \left(
\pi \epsilon \right) \right) +2\frac{1}{\epsilon } \\
&=&-2\left( \frac{1}{\epsilon }\right) +2\frac{1}{\epsilon } \nonumber\\
&=&0 \ , \nonumber\
\end{eqnarray}
so 
\begin{eqnarray}
&&\frac{d}{ds}\left( \left( s-k\right) ^{2}g\left( s\right) \right) \nonumber\\
&&=\frac{\left( -1\right) ^{s}z^{s}\Gamma \left( -s\right) \Gamma \left(
-N+s\right) }{\Gamma \left( 1+n+s\right) }\left( s-k\right)^{2}\nonumber\\
&&\ \ \times\left( {\psi}\left( -N+s\right) -\psi\left( -s\right) -\psi\left(
n+1+s\right) +\ln \left( -z\right) +2\frac{1}{s-k}\right) \nonumber\\
&&=\frac{\left( -z\right) ^{s}\Gamma \left( -s\right) \Gamma \left(
-N+s\right) }{\Gamma \left( 1+n+s\right) }\left( s-k\right) ^{2}\nonumber\\
&&\ \ \times\left({\psi}\left( 1+N-s\right) -\psi\left( 1+s\right) -\psi\left(
1+n+s\right) +\ln \left( -z\right) \right) \ .
\end{eqnarray}

The overlapping pole singularities of the two gamma factors follow from 
\begin{eqnarray}
\Gamma \left( -s\right)  &=&-\frac{1}{\Gamma \left( 1+s\right) }\frac{\pi }{%
\sin \left( \pi s\right) } \nonumber\\
&\approx &\left( -1\right) ^{k+1}\frac{1}{k!}\frac{1}{\epsilon } \\
\Gamma \left( -N+s\right)  &=&-\frac{1}{\Gamma \left( 1+N-s\right) }\frac{%
\pi }{\sin \left( \pi \left( N-s\right) \right) } \nonumber\\
&\approx &\left( -1\right) ^{N-k}\frac{1}{\left( N-k\right) !}\frac{1}{%
\epsilon } \ .
\end{eqnarray}
Overall, this gives residues, for the double poles, with for $k=0,1,\cdots ,N,$%
\begin{eqnarray}
&&\lim_{s\rightarrow k}\frac{d}{ds}\left( \left( s-k\right) ^{2}g\left(
s\right) \right) \nonumber \\
&&=\frac{\left( -z\right) ^{k}}{\left( n+k\right) !}\frac{\left( -1\right)
^{N+1}}{k!\left( N-k\right) !}\left( \psi\left( 1+N-k\right) -{\psi}\left( 1+k\right) -\psi\left( 1+n+k\right) +\ln \left| z\right|
\right) \ . \nonumber\\
\end{eqnarray}
The single-pole residues for $s=k+\epsilon $ for small $\epsilon $ and $%
k=N+1,2,\cdots ,\infty $ come from the poles in $\Gamma \left( -s\right) :$
\begin{eqnarray}
&&\lim_{s\rightarrow k}\left( \left( s-k\right) g\left( s\right) \right)   \nonumber\\
&&=\lim_{s\rightarrow k}\left( \left( s-k\right) \frac{\Gamma \left(
-N+s\right) \Gamma \left( -s\right) }{\Gamma \left( 1+n+s\right) }\left(
-z\right) ^{s}\right)   \nonumber\\
&&=\lim_{s\rightarrow k}\left( \left( s-k\right) \frac{\Gamma \left(
-N+s\right) }{\Gamma \left( 1+n+s\right) }\frac{\left( -1\right) ^{k+1}}{%
k!\left( s-k\right) }\left( -z\right) ^{s}\right)   \nonumber\\
&&=-\frac{\left( k-N-1\right) !}{\left( n+k\right) !}\frac{1}{k!}z^{k} \ .
\end{eqnarray}

As a result, the Cauchy integral satisfies 
\begin{eqnarray}
\frac{1}{2\pi i}&&\hspace{-6mm}\oint_{C^{^{\prime }}}\frac{\Gamma \left( a+s\right)
\Gamma \left( -s\right) }{\Gamma \left( b+s\right) }\left( -z\right) ^{s}ds \nonumber\\
&&=\sum_{k=1}^{n}\lim_{s\rightarrow -k}\left[ \left( s+k\right) g\left(
s\right) \right] +\sum_{k=0}^{N}\lim_{s\rightarrow k}\frac{d}{ds}\left(
\left( s-k\right) ^{2}g\left( s\right) \right)\nonumber\\
&&\hspace{16mm} +\sum_{k=N+1}^{\infty}\lim_{s\rightarrow k}\left[ \left( s-k\right) g\left( s\right) \right]   \nonumber\\
&&=\left( -1\right) ^{N}\sum_{k=1}^{n}\frac{\left( k-1\right) !}{\left(
N+k\right) !\left( n-k\right) !}z^{-k}  \nonumber\\
&&\ \ +\left( -1\right) ^{N+1}\sum_{k=0}^{N}\frac{\left( -z\right) ^{k}}{\left(
n+k\right) !}\frac{1}{k!\left( N-k\right) !}\nonumber \\
&&\ \ \ \ \ \ \ \ \times\left( \psi\left(
1+N-k\right) -\psi\left( 1+k\right) -\psi\left( 1+n+k\right)
+\ln \left| z\right| \right)  \nonumber \\
&&\ \ -\sum_{k=N+1}^{\infty }\frac{\left( k-N-1\right) !}{\left( n+k\right) !}%
\frac{1}{k!}z^{k} \ .\label{series-sol}
\end{eqnarray}

We now compare this series solution \eqref{series-sol} with the expression in DLMF \cite[13.2.28]{DLMF} given by
\begin{eqnarray}
\Psi _{DL}&&\hspace{-8mm}\left( N,n,z\right)  \nonumber \\
&=&\sum_{k=1}^{n}\frac{n!\left( k-1\right) !}{\left( n-k\right) !\left(
1-a\right) _{k}}z^{-k}  \nonumber\\
&&-\sum_{k=0}^{-a}\frac{\left( a\right) _{k}}{\left( n+1\right) _{k}k!}%
z^{k}\left( \ln z+\psi \left( 1-a-k\right) -\psi \left( 1+k\right) -\psi
\left( 1+n+k\right) \right)   \nonumber\\
&&+\left( -1\right) ^{1-a}\left( -a\right) !\sum_{k=1-a}^{\infty }\frac{%
\left( k-1+a\right) !}{\left( n+1\right) _{k}k!}z^{k} \ .
\end{eqnarray}
which is 
\begin{eqnarray}
\Psi _{DL}&&\hspace{-10mm}\left( N,n,z\right)   \nonumber\\
&=&\sum_{k=1}^{n}\frac{n!\left( k-1\right) !}{\left( n-k\right) !\left(
1+N\right) _{k}}z^{-k}  \nonumber\\
&&-\sum_{k=0}^{N}\frac{\left( -N\right) _{k}}{\left( n+1\right) _{k}k!}%
z^{k}\left( \ln z+\psi \left( 1+N-k\right) -\psi \left( 1+k\right) -\psi
\left( 1+n+k\right) \right)   \nonumber\\
&&+\left( -1\right) ^{1+N}N!\sum_{k=1-a}^{\infty }\frac{\left( k-N-1\right) !%
}{\left( n+1\right) _{k}k!}z^{k} \label{psi-dl}
\end{eqnarray}
or 
\begin{eqnarray}
\Psi _{DL}&&\hspace{-10mm}\left( N,n,z\right)  \nonumber \\
&=&N!n!\sum_{k=1}^{n}\frac{\left( k-1\right) !}{\left( n-k\right) !\left(
N+k\right) !}z^{-k}  \nonumber\\
&&-N!n!\sum_{k=0}^{N}\frac{\left( -1\right) ^{k}}{\left( N-k\right) !\left(
n+k\right) !k!}z^{k} \nonumber \\
&&\ \ \ \ \ \ \ \ \times\left( \ln z+\psi \left( 1+N-k\right) -\psi \left(
1+k\right) -\psi \left( 1+n+k\right) \right)  \nonumber \\
&&+\left( -1\right) ^{1+N}N!n!\sum_{k=1-a}^{\infty }\frac{\left(
k-N-1\right) !}{\left( n+k\right) !k!}z^{k} \ .
\end{eqnarray}
Comparing the DLMF series with our contour integral result \eqref{series-sol}, we have 
\begin{equation}
\Psi _{DL}\left( N,n,z\right)
=\left( -1\right) ^{N}N!n!\frac{1}{2\pi i}\oint_{C^{^{\prime }}}\frac{%
\Gamma \left( a+s\right) \Gamma \left( -s\right) }{\Gamma \left( b+s\right) }%
\left( -z\right) ^{s}ds 
\end{equation}
or
\begin{eqnarray}
\Psi _{DL}\left( N,n,z\right)  &=&\frac{\Gamma \left( 1+n\right) }{%
\epsilon\Gamma \left( -N+\epsilon \right) }\frac{1}{2\pi i}\oint_{C^{^{\prime }}}%
\frac{\Gamma \left( a+s\right) \Gamma \left( -s\right) }{\Gamma \left(
b+s\right) }\left( -z\right) ^{s}ds  \nonumber\\
&=&\frac{\Gamma \left( b\right) }{\epsilon \Gamma \left( a+\epsilon \right) }%
\frac{1}{2\pi i}\oint_{C^{^{\prime }}}\frac{\Gamma \left( a+\epsilon
+s\right) \Gamma \left( -s\right) }{\Gamma \left( b+s\right) }\left(
-z\right) ^{s}ds  \nonumber\\
&=&\lim_{\epsilon\rightarrow 0}\frac{1}{\epsilon}\frac{1}{2\pi i}%
\oint_{C^{^{\prime }}}\frac{\Gamma \left( a+s\right) }{\Gamma \left(
a\right) }\frac{\Gamma \left( b\right) }{\Gamma \left( b+s\right) }\Gamma
\left( -s\right) \left( -z\right) ^{s}ds  \nonumber\\
&=&\lim_{\epsilon \rightarrow 0}\frac{1}{\epsilon} \frac{1}{2\pi i}%
\oint_{C^{^{\prime }}}\frac{\left( a\right) _{s}}{\left( b\right) _{s}}%
\Gamma \left( -s\right) \left( -z\right) ^{s}ds \ .\nonumber \\
\end{eqnarray}
wherein $a=-N+\epsilon$.

The above constitutes a derivation and verification of the DLMF series solution given in \cite[13.2.28]{DLMF}.

\section{Comparison of our closed-form solution with the the series solution}
\label{compare}

Our closed-form second solution \eqref{psibar}, for $a=-N,\ \ b=n+1$, is 
\begin{equation}
\Psi _{PM}(N,n,x)=\frac{n!}{(N+n)!}\,P(N,n,x)\frac{e^{x}}{x^{n}}+\Phi \left(
N,n,x\right) \int_{-x}^{\infty }\frac{e^{-s}}{s}\,ds \ .
\end{equation}
where
\begin{eqnarray}
&&P\left( N,n,x\right)  \nonumber\\
&&=\sum_{m=0}^{n-1}\frac{\left( N+m\right) !\left(
n-m-1\right) !}{m!}x^{m} \nonumber\\
&&\ \ -\,x^{n}\sum_{m=0}^{N-1}\sum_{k=0}^{N-m-1}\frac{N!}{%
\left( N-k-m-1\right) !}\frac{(N+n)!}{\left( n+k+m+1\right) !}\frac{%
(-1)^{k}k!}{\left( k+m+1\right) !}\left( -x\right) ^{m}. \nonumber\\
\end{eqnarray}
and
\begin{equation}
\Phi (N,n,x)=\sum_{k=0}^{N}\frac{N!n!}{\left( N-k\right) !(n+k)!k!}\left(
-x\right) ^{k} \  ,
\end{equation}
which is the first (regular) solution $\,_1F_1(-N;n+1;x)$.

One can see that $\Psi _{PM}(N,n,x)$ coincides with $\Psi _{DL}\left(
N,n,x\right) $ as follows. If one expands $\Psi _{PM}(N,n,x)$ in a Laurent
power series (with possible logarithmic terms), the only terms that have
negative powers in $x$  come from
\begin{equation*}
\frac{n!}{(N+n)!}\frac{e^{x}}{x^{n}}\sum_{m=0}^{n-1}\frac{\left( N+m\right)
!\left( n-m-1\right) !}{m!}x^{m}
\end{equation*}
Expanding, but keeping only powers up to $x^{n-1},$ we will have 
\begin{eqnarray*}
&&\frac{n!}{(N+n)!}e^{x}\sum_{m=0}^{n-1}\left[ \frac{\left( N+m\right) !\left(
n-m-1\right) !}{m!}\right] x^{m}\\
&=&\frac{n!}{(N+n)!}\sum_{m=0}^{n-1}%
\sum_{k=0}^{n-1-m}\frac{\left( N+m\right) !\left( n-m-1\right) !}{m!k!}%
x^{m+k} \\
&=&\frac{n!}{(N+n)!}\sum_{m=0}^{n-1}\sum_{k=0}^{n-1-m}\frac{\left(
N+m\right) !\left( n-m-1\right) !}{m!k!}x^{m+k} \\
&=&\frac{n!}{(N+n)!}\sum_{l=0}^{n-1}\left( \sum_{k=0}^{l}\frac{\left(
N+l-k\right) !\left( n-1-l+k\right) !}{\left( l-k\right) !k!}\right) x^{l} \\
&=&\frac{n!}{(N+n)!}\sum_{l=0}^{n-1}\left( \frac{\left( N+l\right) !\left(
n-l-1\right) !}{l!}\sum_{k=0}^{l}\frac{\left( -l\right) _{k}\left(
n-l\right) _{k}}{\left( -N-l\right) _{k}k!}\right) x^{l} \\
&=&\frac{n!}{(N+n)!}\sum_{l=0}^{n-1}\left( \frac{N!\left( N+n\right) !\left(
n-1-l\right) !}{\left( N+n-l\right) !l!}\right) x^{l} \\
&=&\sum_{l=0}^{n-1}\frac{n!N!\left( n-1-l\right) !}{\left( N+n-l\right) !l!}%
x^{l}
\end{eqnarray*}
where we used the Gauss identity: $\,_{2}F_{1}\left( -l,b;c;1\right) =\left(
c-b\right) _{l}/\left( c\right) _{l}$. Comparing to the DLMF first term
(times $x^{n}$): 
\begin{equation*}
\sum_{k=1}^{n}\frac{n!(k-1)!}{(n-k)!(1+N)_{k}}x^{n-k}=\sum_{k=1}^{n}\frac{%
N!n!(k-1)!}{(n-k)!(N+k)!}x^{n-k}=\sum_{l=0}^{n-1}\frac{N!n!(n-l-1)!}{%
l!(N+n-l)!}x^{l}\,,
\end{equation*}
we find they agree in sign and magnitude. 

The logarithmic terms in $\Psi _{PM}(N,n,x)$ come from
\begin{equation*}
\Phi \left( N,n,x\right) \int_{-x}^{\infty }\frac{e^{-s}}{s}\,ds \ .
\end{equation*}
Using the expansion
\begin{equation*}
\int_{-x}^{\infty }\frac{e^{-s}}{s}\,ds={Ei}\left( 1,-x\right) =-\ln
x-\gamma -\sum_{k=1}^{\infty }\frac{1}{k\,k!}x^{k}
\end{equation*}
and comparing with the DLMF logarithmic term (see \eqref{psi-dl})
\begin{equation*}
-\Phi \left( N,n,z\right) \left( \ln z\right) \ ,
\end{equation*}
we see that the log term in $\Psi _{PM}(N,n,x)$ matches that in $\Psi
_{DL}(N,n,x)$. Moreover, the $\gamma $ term, coming from the $\psi $ terms
in $\Psi _{DL}(N,n,x),$ also matches.

As $\Psi _{PM}(N,n,x)$ and $\Psi _{DL}(N,n,x)$ are both solutions to the
confluent hypergeometric equation, they can only differ by another
independent solution, i.e.
\begin{equation}
\Psi _{PM}(N,n,x)=c_{1}\Psi _{DL}(N,n,x)+c_{2}\Phi \left( N,n,x\right) 
\end{equation}
By finding that $\Psi _{PM}(N,n,x)$ matches $\Psi _{LM}(N,n,x)$ over a range
of $x,$ we have $c_{1}=1$ and $c_{2}$ $=0 \ ,$
establishing that
\begin{equation}
\Psi _{PM}(N,n,x)=\Psi _{DL}(N,n,x) \ . 
\end{equation}

\newpage




\label{lastpage}

\end{document}